\documentclass[11pt]{articulofederico}
\usepackage{graphicx}
\usepackage{dcolumn}
\usepackage{enumerate}
\usepackage{subfigure}
\usepackage{epsfig}
\usepackage{amssymb}
\usepackage{epsf} 	
\usepackage{amscd}

\usepackage[latin1]{inputenc}
\usepackage[spanish, activeacute]{babel}

\usepackage{xcolor} % Para los colores

\usepackage{amsfonts} 

\usepackage{amssymb}
\usepackage{enumerate}
\usepackage{graphicx}

\usepackage{amsmath,amssymb,enumerate,amsthm,flafter}

\newtheorem{thm}{Teorema}[section]
\newtheorem{Teorema}{Teorema}[section]

\newtheorem{obs}[thm]{Observaci\'on} 
\theoremstyle{definition}
 \newtheorem{defn}{Definici\'on}[section]
  \newtheorem{definicion}{Definici\'on}[section]
    
        \newtheorem{Observaci'on}{Observaci\'on}[section]
 
  \newtheorem{Ejemplo}[defn]{Ejemplo}

    \newtheorem{Ejercicio}[defn]{Ejercicio}

\renewcommand\figurename{Figura}

\DeclareMathOperator{\tr}{tr}
\DeclareMathOperator{\sgn}{sgn}

\def\Sym{{Sym}}

\def\S{\mathbb{S}}

\def\C{\mathbb{C}}
\def\Q{\mathbb{Q}}

\def\Z{\mathbb{Z}}

\def\dom{\unrhd}

\def\ch{{\bf ch}}

\def\1{{\bf 1}}

\def\P{{\bf \P}}
\def\H{\mathcal{H}}
\def\L{\mathcal{L}}
\def\gr{R(\S)}
\def\kk{{\mathbb{K}}}
\def\grn{R(\S_n)}
\def\kr{\gamma^{\lambda}_{\mu,\nu}}

\def\Young#1{\vbox{\smallskip\offinterlineskip 
    \halign{&\vbox{##}\kern-\Thickness\cr #1}}} 
 
\newdimen\Squaresize \Squaresize=14pt 
\newdimen\Thickness \Thickness=.3pt 
\newdimen\Correction \Correction=7pt 
 
\def\Vide#1{\hbox{ 
       \vbox to \Squaresize{\vss 
          \hbox to \Squaresize{\hss#1 \hss}\vss} 
    \hskip-\Correction} 
   \kern-\Thickness}

\def\Carre#1{\hbox{\vrule width \Thickness 
   \vbox to \Squaresize{\hrule height \Thickness\vss 
      \hbox to \Squaresize{\hss#1\hss} 
   \vss\hrule height\Thickness} 
   \unskip\vrule width \Thickness} 
   \kern-\Thickness}

\def\Box#1{\Carre{$\scriptstyle#1$}}

\renewcommand{\refname}{\textsf{Referencias}}
\renewcommand{\abstractname}{\textsf{Resumen}}

\begin{document}
\title{\textsf{Tres lecciones en combinatoria algebraica.  \\
\normalsize{II. Las funciones sim\'etricas y la teor'ia de representaciones.}}}
%\title{\textsf{Lifted generalized permutahedra and \\ composition polynomials.}}

 \author{\textsf{Federico Ardila\footnote{\textsf{San Francisco State University, San Francisco, CA, USA y Universidad de Los Andes, Bogot\'a, Colombia, federico@sfsu.edu -- financiado por la CAREER Award DMS-0956178 y la beca DMS-0801075 de la National Science Foundation de los Estados Unidos, y por la SFSU-Colombia Combinatorics Initiative.}}}\qquad
\textsf{Emerson Le\'on\footnote{\textsf{Freie Universit\"at Berlin, Alemania,  emerson@zedat.fu-berlin.de -financiado por el Berlin Mathematical School.}}}\\
\textsf{Mercedes Rosas\footnote{\textsf{Universidad de Sevilla, Espa\~na, mrosas@us.es -- financiada por los proyectos MTM2007--64509 del Ministerio de Ciencias e Innovaci\'on de Espa\~na y FQM333 de la Junta de Andalucia.}}}\qquad
\textsf{Mark Skandera\footnote{\textsf{Lehigh University, %, 14 East Packer Ave, 
Bethlehem, PA, USA, % 18015, 
mas906@math.lehigh.edu -- financiado por la beca H98230-11-1-0192 de la National Security Agency de los Estados Unidos.}}}
}
\date{}
\maketitle

\begin{abstract} % Short abstract
En esta serie de tres art\'\i culos, damos una exposici\'on de varios resultados y problemas abiertos en tres \'areas de la combinatoria algebraica y geom\'etrica:  las matrices totalmente no negativas,
las representaciones
del grupo sim\'etrico $S_n$, y los arreglos de hiperplanos. Esta segunda parte trata la conecci\'on entre las funciones sim\'etricas y la teor\'{\i}a de representaciones.
\end{abstract}

%\tableofcontents

%%%%%%%%%%%%%%%%%%%%%%%%%%%%%%%%%%%%%%%%%%%%%%%%%
%%%%%%%%%%%%%%%%%%%%%%%%%%%%%%%%%%%%%%%%%%%%%%%%%
%%%%%%%%%%%%%%%%%%%%%%%%%%%%%%%%%%%%%%%%%%%%%%%%%
%%%%%%%%%%%%%%%%%%%%%%%%%%%%%%%%%%%%%%%%%%%%%%%%%

En marzo de 2003 se llev\'o a cabo el Primer Encuentro Colombiano de Combinatoria en Bogot\'a, Colombia. Como parte del encuentro, se organizaron tres minicursos, dictados por Federico Ardila, Mercedes Rosas, y Mark Skandera. Esta serie resume el material presentado en estos cursos en tres art\'{\i}culos:   \emph{I. Matrices
totalmente no negativas y funciones sim\'etricas} \cite{ALRS1}, 
\emph{II. Las funciones sim\'etricas y la teor\'{\i}a de las representaciones} \cite{ALRS2}, y \emph{III. Arreglos de hiperplanos.} \cite{ALRS3} 

%
%En la primavera del 2003, la Universidad de los Andes (Bogot'a, Colombia) nos invit'o a dar un ciclo de conferencias sobre
%la combinatoria algebraica. Fue una ocasi'on inolvidable, tanto por el placer de trabajar con sus estudiantes, excepcionalmenete brillantes y
%motivados, como por la belleza de su campus y la hospitalidad de su gente. 

En esta segunda lecci'on hablaremos sobre nuestra
introducci'on a la teor'ia de representaciones del grupo sim'etrico, del grupo lineal general,  y a la teor'ia de las
 funciones sim'etricas. Esperamos que este trabajo haga justicia a la experiencia
vivida en esos d'ias.

\section{\textsf{Algunas nociones b'asicas.}}

%Con la finalidad de que esta lecci'on se pueda leer de manera independiente de nuestra primera lecci'on, volvemos a introducir ahora algunas lecciones
%centrales de la combinatoria algebraica

Empezamos por introducir algunas nociones b'asicas de combinatoria algebraica. 
%Algunas de ellas ya aparecen en la primera lecci'on
%pero volveremos a hablar de ellas para que sea posible leerla de manera independiente.
Una sucesi\'on debilmente decreciente  de enteros no negativos
  se  denomina {\em partici\'on} y se escribe como
$ \lambda = (\lambda_1, \lambda_2, \dotsc).$ Un ejemplo de partici'on es la siguiente $\lambda=(7,7,4,1,1,1,0,0)$
Aquellos $\lambda_i$ mayores que cero se llaman las {\em partes} de
$\lambda$, y decimos que dos particiones son iguales si difieren solamente en el n'umero de ceros.

El n'umero de partes de $\lambda$ se denomina la longitud de $\lambda$ y se denota por $\ell(\lambda)$.
En ocasiones escribimos la partici'on $\lambda$ describiendo la multiplicidad de sus partes. Por ejemplo,
para nuestro ejemplo, escribimos
$\lambda=(1^3\, 4\, 7^2).$ 

Decimos que $\lambda$ es una {\em partici\'on de $n$} si  $
\lambda_1 + \cdots + \lambda_{\ell(\lambda)} = n,.
$ En este caso, escribimos $\lambda \vdash n$, o $|\lambda| = n$.
Identificamos una partici'on de $n$ con su {\em diagrama de Young}, un arreglo de filas de cuadrados, justificados por la derecha,
donde la $i$-'esima fila contiene $\lambda_i$ cuadrados, tal y como se ilustra en la  Figura 1. La
{\em partici'on transpuesta} de $\lambda$ se define a trav'es de su diagrama de Young, que se
obtiene al reflejar el diagrama de $\lambda$ sobre su diagonal principal. Se denota por $\lambda'$.
\begin{figure}[h]
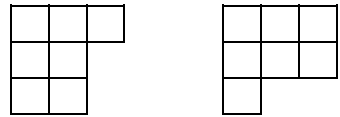

\begin{eqnarray*}
&\begin{matrix} \Young{ 
	\Box{}& \Box{}& \Box{} \cr  
	\Box{}& \Box{} \cr 
		\Box{}& \Box{} \cr 
 }
     \end{matrix}  
     \,\,\,\,\,\,\,\,\,\,\,\,\,\,\,
     &\begin{matrix} \Young{ 
	\Box{}& \Box{} & \Box{}\cr  
	\Box{}& \Box{}  & \Box{}\cr 
	\Box{}\cr }
     \end{matrix} 
      \end{eqnarray*}
\caption{   Los diagramas de Young de la partici'on $\lambda=(3,2,2)$ y de su transpuesta $\lambda'=(3,3,1)$.}
\end{figure}

Similarmente, dadas dos particiones $\mu$ y $\lambda$, cuyos diagramas de Young satisfacen que $\mu \subseteq \lambda$, definimos el tablero de Young sesgado, $\lambda / \mu$
como aquel que se obtiene al restar al diagrama de $\lambda$ el diagrama de $\mu$ (vistos como conjuntos). 
\begin{figure}[h]
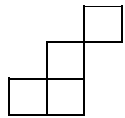

\begin{eqnarray*}
&\begin{matrix} \Young{ 
	& & \Box{} \cr  
	& \Box{} \cr 
		\Box{}& \Box{} \cr 
 }
     \end{matrix}  
       \end{eqnarray*}
\caption{ El diagrama de Young  sesgado correspondiente a  $(3,2,2) / (2,1)$.}
\end{figure}

Sea $\lambda$ una partici'on de  $n$, un {\em tableau} es una manera de asignar un n'umero en $\mathbb{N}$ a cada celda del diagrama de $\lambda$, 
donde  es posible utilizar el mismo número repetidamente. En algunas ocasiones, es conveniente pedir que todas las entradas utilizadas 
pertenezcan a $[n]$ para alg'un $n$, donde $[n]$ denota al conjunto $\{1,2,\cdots,n\}$.
%Esto equivale a trabajar con un conjunto finito de variables en el caso de las funciones sim'etricas que 
%estudiaremos posteriormente.

Un  tableau es {\em semi-estándar} si sus columnas crecen estrictamente, mientras que sus filas crecen débilmente. La sucesión \[
(\alpha_1(T), \alpha_2(T), \cdots, \alpha_{n_{\ell}}(T))
\]
 donde $\alpha_i(T)$  es número de veces que aparece el número $i$ en el tableau $T$ se denomina el {\em contenido} de $T$.
Un tableau semi-est'andar es {\em estándar} si cada uno de los números en $[n]$ aparece exactamente una vez.

%En algunos libros, como por ejemplo \cite{F}, se refieren a los tableaux semi-est'andar
%simplemente como tableaux.

Sea $T$ un tableau de forma $\lambda$. El {\em peso} de $T$, que denotamos por $x^T$ se define como
\[
x^{T}=x_1^{\alpha_1(T)}x_2^{\alpha_2(T)}\ldots
\]

El {\em coeficiente de Kostka} $K_{\lambda, \mu}$ se define como el número de tableaux semi-estándar de forma $\lambda$ y contenido $\mu$.
En particular, denotamos por $f^{\lambda}$ al n'umero de tableaux est'andar de forma $\lambda$.    Esto es,  $f^{\lambda}=K_{\lambda,(1^n)}$.

\begin{figure}[h]\label{fig.semis}
\begin{eqnarray*}
&\begin{matrix} \Young{ 
	\Box{1}& \Box{1}& \Box{2} \cr  
	\Box{2}& \Box{3} \cr }
     \end{matrix} 
&\,\,\,\,\,\,\,\,\,\,\,\,\,\,\,\,\,\,\,\,\,\,\,\,\,\,\,\,\,\,
\begin{matrix} \Young{ 
	\Box{1}& \Box{1} & \Box{3}\cr  
	\Box{2}& \Box{2} \cr }
     \end{matrix} 
 \end{eqnarray*}
\caption{   Los dos tableaux semi-estándar de forma $(3,2)$ y contenido $(2,2,1)$. Ambos tienen peso $x_1^2x_2^2x_3$.
}
\label{semistandar}
\end{figure}

%Existe una f'ormula cerrada que nos permite calcular $f^{\lambda}$, llamada en ingl'es de hook-f'ormulas (f'ormula de las escuadras). Por otra parte,
%para calcular el n'umero de tableaux semi-est'andar de forma $\lambda$ y entradas en $[n]$ tenemos una funci'on generatriz, llamada en 
%ingl'es la hook-content, descrita por Stanley. El n'umero de tableaux semi-est'andar de forma $\lambda$ tambi'en puede expresarse como
% una especializaci'on de las funciones de Schur $s_{\lambda}$, ver \cite{Au, EC1}.

Sean $\mu$ y $\nu$ dos particiones del mismo entero. Definimos al {\em orden de dominancia}, denotados por $ \dom$, diciendo que $\mu  \dom \nu$ si para 
cada $k$ tenemos que 
\[
\mu_1+\mu_2+\cdots+\mu_k \ge \lambda_1+\lambda_2+\cdots+\lambda_k
\]
En f'isica, este orden se conoce como {\em mayorizaci'on}.

\begin{Ejercicio}[Los coeficientes de Kostka y el orden de dominancia]
Supóngase que $\mu$ y $\lambda$ son particiones de un mismo entero y que $K_{\lambda,\mu}  \neq 0$. Demuéstrese que  $\lambda \dom \mu$ en el orden de dominancia. Demuéstrese, además, que $K_{\lambda, \lambda}=1$.
\end{Ejercicio}

Trabajaremos siempre sobre el cuerpo de los n'umeros complejos.
Recordemos que el {\em grupo simétrico ${\S}_n$} se define como el conjunto de todas las permutaciones del conjunto $[n]$, junto con la operación de composición.  A los elementos del grupo simétrico  los denominamos permutaciones. Denotamos por $\pi \sigma$ la permutación obtenida al aplicar primero $\sigma$ y luego $\pi$ a los elementos de $[n]$.  A la permutación identidad la denotamos con la letra $\epsilon$.

Utilizamos dos maneras diferentes para denotar permutaciones. De acuerdo con la primera convención, la permutación $\pi$ que  envía el número $i$ en $\pi(i)$ la denotamos con la  palabra $\pi(1) \, \pi(2) \, \cdots  \pi{(n)} $. De acuerdo a la segunda, escribimos $\pi$ como el producto de sus ciclos.   
Por ejemplo, la permutación  $\pi$  definida por  $\pi(1)=2,$ $\pi(2)=3,$ $\pi(3)=7,$ $\pi(4)=4,$ $\pi(5)=1,$ $\pi(6)=8,$ $\pi(7)=5,$ $\pi(8)=6$ se escribe, utilizando la primera convención, por la palabra  $2 \, 3 \,7 \,4 \,1 \, 8 \, 5 \, 6$. Su descomposición en ciclos es:  $(1\,  2 \, 3 \, 7 \, 5 ) (4) (6 \, 8)$. 
El {\em tipo} de una permutación  es la partición definida por las longitudes de los ciclos que aparecen en ella. La   permutación $\pi$ que aparece en nuestro ejemplo tiene tipo $5 2 1 \vdash 8$.  

%A lo largo de est'a lecci'on asignaremos al lector una breve lista de ejercicios.
A cada partici'on $\lambda$ le asociamos el entero $z_{\lambda} $ definido como 
\[
z_{\lambda} = 1^{m_1} m_1!\, 2^{m_2} m_2! \, \cdots n^{m_n} m_n!.
\]

 \begin{Ejercicio} \label{index}  Demostrar que dos permutaciones son conjugadas si y s'olo si tienen el mismo tipo. M'as a'un, demostrar que el número de permutaciones de ${\S}_n$ cuyo  tipo es la partición $\lambda=1^{m_1}\, 2^{m_2} \, \cdots n^{m_n}$ de $n$ es igual a 
$
 \frac{n!}{z_{\lambda}}.
 $
%Similarmente, demostrar que $z_{\lambda}$ es igual al n'umero de permutaciones en $\mathbb{S}_n$ que conmutan con una permutaci\'on fija de tipo $\lambda$.
\end{Ejercicio}

% 
% Si $\pi= \pi_1 \pi_2 \cdots \pi_k $ es una expresión de $\pi$ como producto de transposiciones, definimos  el {\em signo} de $\pi$ como $\sgn(\pi)=(-1)^k$. El signo de una permutación no depende de su descomposición como producto de transposiciones. Nótese que $\sgn(\pi) \sgn(\tau) = \sgn(\pi \tau)$ y que $\sgn(\epsilon) =1.$
 
 %%%%%%%%%%%%%%%%%%%%%%%%%%%%%%%%%%%%%%%%%%%%%%%%%
%%%%%%%%%%%%%%%%%%%%%%%%%%%%%%%%%%%%%%%%%%%%%%%%%
%%%%%%%%%%%%%%%%%%%%%%%%%%%%%%%%%%%%%%%%%%%%%%%%%
%%%%%%%%%%%%%%%%%%%%%%%%%%%%%%%%%%%%%%%%%%%%%%%%%

\section{\textsf{Introducci'on a la teor'ia de representaciones de grupos.}}
Comenzamos con una breve exposición de la teoría general de representaciones de grupos, centrándonos en el caso del grupo simétrico y basada en los trabajos de  Fran\c cois Bergeron \cite{B}, William Fulton \cite{F} y Bruce Sagan \cite{Sagan}.  %Nuestro objetivo es que el lector observe el paralelismo que existe entre esta
%teor'ia de representaciones y las funciones sim'etricas. 

Sea $G$ un grupo, una {\em representación (matricial, compleja) de $G$}  es un homomorfismo de grupos  entre  $G$ y el grupo de las matrices invertibles de orden dado,   $GL_d=GL_d(\mathbb C)$. 
Esto es,
\begin{align*}
X : G &\rightarrow GL_d\\
 \pi &\mapsto X(\pi).
\end{align*}
Al parámetro $d$ lo llamamos {\em orden o dimensi'on} de la representación. Dada una representación matricial $X$,  denotamos $X(\pi) \vec v$ por $\pi \cdot \vec v$.
En esta situación decimos que  $G$ {\em actúa (linealmente)} sobre $\C^d$.

Consideramos ahora dos acciones del grupo sim'etrico sobre $[m]^n$:
 \begin{align*}
\sigma \cdot (v_1, v_2, \cdots, v_n) & =  (\sigma(v_1), \sigma(v_2), \cdots, \sigma(v_n)), &&\text{$\sigma \in \S_m$ permuta las entradas de $v$.}\\
 (v_1, v_2, \cdots, v_n) \cdot \tau &=  (v_{\tau(1)},v_{\tau(2)}, \cdots, v_{\tau(n)}), &&\text{$\tau \in \S_n$ permuta las posiciones en $v$.}
\end{align*}

\begin{Ejercicio} Demostrar que las dos definiciones que acabamos de dar son, efectivamente, acciones de $\S_m$ (respectivamente $\S_n$)
sobre $[m]^n$. Por otra parte, verificar que si definimos
$\tau \cdot (v_1, v_2, \cdots, v_n)  =  (v_{\tau(1)},v_{\tau(2)}, \cdots, v_{\tau(n)}), $ no obtenemos una acci'on de $\S_n$ sobre
$[m]^n$ ya que no proviene de un homomorfismo de grupos.
\end{Ejercicio}

Dos ejemplos de representaciones del grupo simétrico son la {\em representación trivial}, que se obtiene al enviar todos los elementos de 
 $\S_n$ a la matriz identidad, y la {\em representación  alternante} definida por $X(\pi)= (\sgn (\pi))$.  Ambas representaciones tienen orden uno.

Otras dos representaciones del grupo simétrico particularmente importantes  son la {\em representación definición}, obtenida al hacer actuar $\S_n$ sobre el conjunto $[n]$ de la manera canónica (utilizando como base a los elementos de $[n]$), y la {\em representación regular}, que proviene de hacer actuar un grupo finito $G$ sobre sí mismo (utilizando como base a los elementos de $G$).  

Un {\em homomorfismo de representaciones} es una transformaci'on lineal $Y : \mathbb{C}^m \to \mathbb{C}^d$ tal que $Y(g v) = g Y(v)$ para todo $g \in G$ y $v \in \mathbb{C}^m$. 
Si $Y$ es invertible tenemos entonces un  {\em isomorfismo.}

\begin{Ejercicio} [La representación  definición de ${\S}_3$.]  \label{matricesS3}
Demuestre que si hacemos actuar ${\S}_3$ sobre  $\C^3$ permutando los vectores de la base canónica obtenemos  la siguiente  representación matricial 
\[
\begin{array}{llcr}  
&X( 1\, 2\, 3)= 
\begin{pmatrix}
 1 & 0 & 0 \\ 
 0 & 1 & 0 \\
 0 & 0 & 1
 \end{pmatrix}
 &X( 2\, 1\, 3 )=
 \begin{pmatrix}
 0 & 1 & 0 \\ 
 1 & 0 & 0\\
 0 & 0 & 1
 \end{pmatrix}
 &X(3\, 2\, 1)=
 \begin{pmatrix}
 0&0&1\\
 0&1&0\\
 1&0&0
 \end{pmatrix}    
 \\ %%nueva linea 
 &X( 1\, 3\, 2 )=
 \begin{pmatrix}
 1&0&0\\
 0&0&1\\
 0&1&0
 \end{pmatrix}
 &X(2\, 3\, 1 )=
 \begin{pmatrix}
 0&0&1\\
 1&0&0\\
 0&1&0
 \end{pmatrix}
 &X( 3\, 1\, 2 )=
 \begin{pmatrix}
 0&1&0\\
 0&0&1\\
 1&0&0
 \end{pmatrix}
\end{array}
\]
 Calcule las matrices correspondientes a la representación definici'on de $\S_4$ sobre  $\C^4$,
 y observe  que, para cada $n$ las matrices obtenidas mediante esta construcci'on son matrices ortogonales, y en
  consecuencia sus inversas vienen dadas por sus transpuestas.
 
\end{Ejercicio}

En general, decimos que una representación de un grupo finito $G$ es una {\em representación por permutaciones}  cuando proviene de la acción de 
$G$ sobre un conjunto finito permutando sus elementos.  Esta acci'on nos proporciona
un homomorfismo natural de  $G$ al grupo de las matrices de permutaci'on,  un interesante subgrupo de $GL_d$.
Las matrices correspondientes a representaciones por permutaciones son matrices de permutaciones (matrices en las que cada fila y cada columna contiene exactamente una entrada diferente de cero e igual a uno).
Tanto la representación definición como la representación regular son ejemplos de representaciones por 
permutaciones. Por otra parte, la representación alternante no lo es. 

Las representaciones por permutaciones son particularmente interesantes desde
un punto de vista combinatorio. Por ejemplo, la traza de $ X(\pi)$ cuenta el n'umero de elementos del conjunto que permanecen fijos bajo la acci'on de $X$. En
particular,  el orden (o dimensi'on) de una representaci\'on  por permutaciones $X$ viene dado por la traza de $X(\epsilon)$.
%Para introducir una familia de objetos, es necesario definir la naturaleza de las aplicaciones entre sus objetos. Para esto introducimos las nociones de ${\S}_n$-homomorfismo y de ${\S}_n$-isomorfismo.

%\color{red} Redacción un poco rara, lo de la familia.
%\color{black}

%\todo{aqui??}

%\begin{definicion}[${\S}_n$--homomorfismo]

%Sean $X$ y $Y$ representaciones de ${\S}_n$ de orden $d$, un {\em ${\S}_n$-homomorfismo} es una  transformación lineal
%\[
%T : {\C}^d \rightarrow {\C}^d
%\]
%tal que $T X(\pi) = Y(\pi) T$, para todo $\pi \in {\S}_n$. Si la transformación lineal $T$ es invertible, decimos entonces que $T$ es  un {\em $ {\S}_n$-isomorfismo.}
%Si existe un ${\S}_n$-isomorfismo entre dos representaciones decimos entonces que son {\em equivalentes.}

%\end{definicion}

%  tExiste un sencillo invariante sobre el conjunto de representaciones de todo grupo finito, que conocemos de nuestro primer curso de \'algebra lineal.

\subsection{\textsf{Los caracteres}}

El ejemplo con el que concluimos la secci'on anterior nos sugiere la importancia de considerar la traza de una representaci'on, y nos pone en contacto con un invariante fundamental de cualquier representaci'on de order finito de un grupo.

\begin{definicion}[Car\'acter] Sea $X$ una representación  matricial de un grupo $X : G \to GL_d(\mathbb{C})$.
El {\em car\'acter} de $X$ es la función $\chi$:
\begin{align*}
\chi : G &\to {\C}\\
 \pi &\mapsto \tr X(\pi).
\end{align*}
donde $\tr X$ denota la traza de la matriz $X$. 
%Algunas veces es conveniente   identificar la función $\chi$ con el vector 
%\[
%\chi=(\chi(\pi_1), \chi(\pi_2), \cdots, \chi(\pi_{n!})) \in {\C}^{n!}.
%\]
%(Para esto es necesario establecer un orden total sobre el conjunto de las permutaciones.) 
\end{definicion}

El car\'acter de una representación no depende de  la base utilizada para construir la matriz asociada a la transformación lineal.  Por otra parte,
no es necesario que el grupo sea finito para la definici'on de car'acter.

%La tabla de caracteres de la representaci'on $X$ de $G$ nos proporciona el valor del car'acter $X$ en todos los elementos de $G$.
Por ejemplo,   el car'acter   de la representaci'on definici'on de $\S_3$, se calcula rápidamente a partir
de los resultados obtenidos en el ejercicio \ref{matricesS3}. 
\begin{align*}
&&\chi_{def}(123)=3  &&\chi_{def}(213)=1 &&\chi_{def}(321)=1\\
&&\chi_{def}(132)=1 &&\chi_{def}(231)=0 &&\chi_{def}(312)=0
\end{align*}
Sabemos que la dimensi'on de la representaci'on definici'on viene dada por el valor de su car'acter en la identidad.
Más generalmente, tenemos que   $\chi_{def}(\sigma)$ cuenta el n'umero de puntos fijos de $\sigma$.

Los caracteres juegan un rol fundamental en la teoría de representaciones de los grupos finitos. 
Esto se debe al siguiente resultado que nos dice que, no solamente son
 un invariante dentro de sus clases de isomorfismos, sino que nos permite distinguirlas.
\begin{Teorema}
Dos representaciones de un grupo finito $G$ tienen el mismo car\'acter s\'i y s\'olo s\'i son isomorfas. 
Esto es, si $X$ es una representación con car\'acter $\chi$ y $Y$ es una representación con car\'acter $\phi$,
\[
X \cong Y \text{ si y s\'olo si } \chi(g) = \phi(g) \text{ para cada $g$ en $G$.} 
\]
\end{Teorema}

Pronto volveremos a la noción de car\'acter. Ahora nos planteamos el problema de 
descomponer una representación  $X : G \to GL_d(\mathbb{C})$
como una suma directa de representaciones más sencillas.  Empezamos por estudiar los subespacios de $\mathbb{C}^d$ que permanecen invariantes bajo $G$.

\begin{definicion}[Subespacio Invariante]
 Un subespacio vectorial $W$ de un espacio $V$ es  {\em invariante bajo $G$} si 
 \[
 \pi \cdot W = X(\pi) \, W = \{ X(\pi)\, (w) : w \in W \}  \subseteq W
 \]
 para cada $\pi$ en $G$.
\end{definicion}

 Por ejemplo, el subespacio  generado por $W=\vec 1 + \vec 2 + \vec 3$ es invariante bajo la representación definición de $\S_3$.
 
El complemento ortogonal  de un subespacio $G$-invariante $W$ tambi'en es invariante, siempre y cuando
el producto escalar que lo define sea invariante bajo $G$, es decir $\langle u, v \, \rangle = \langle g \cdot u, g \cdot v \, \rangle $ 
para todo $u, v \in V$ y $g \in G$. Para la representaci'on definici'on de
$\mathbb{S}_n$, el producto escalar can'onico es invariante. (El que hace que los vectores de la base can'onica sean ortonormales.)
Tenemos entonces que el complemento ortogonal al subespacio $W=\vec 1 + \vec 2 + \vec 3$ tambi'en es invariante por
 la representación definición de $\S_3$. Este espacio est'a  generado por $\vec 2 - \vec 1,$ y $\vec 3 - \vec 1$.
  %\footnote{Para cualquier representaci'on de un grupo finito $G$, a partir de cualquier producto escalar $\langle \cdot, \cdot \, \rangle$
 % construimos un producto escalar invariante bajo la acci'on de $G$ utilizando  la siguiente construcci'on:
% \[
% \langle u,v \rangle = \frac{1}{|G|} \sum_{g \in G} \langle  g \cdot u ,  g \cdot v  \rangle.
% \] Al menos si trabajamos sobre un cuerpo de caracter'istica cero.}. 
%  Hemos encontrado una nueva familia de representaciones del grupo simétrico.
%	 
%	 \todo{Incluir ilustracion.}

La representación definición de $\S_n$ actúa de manera trivial sobre el subespacio generado por $1+2+\ldots+n$. La restricción  de esta representación al complemento ortogonal de este subespacio lo denominamos la {\em representaci'on est'andar. } Tiene por base  $\vec 2 - \vec 1, \vec 3 - \vec 1, \cdots, \vec n - \vec 1$, y en
consecuencia orden $n-1$.

\begin{Ejercicio}[La representaci'on est'andar de $\S_3$] \label{descirre}
Demostrar que si escribimos la transformación lineal $X$ definida en el Ejercicio \ref{matricesS3} en la base
$ \{ \vec 1 + \vec 2 + \vec 3, \, \vec 2 - \vec 1, \, 
 \vec 3 -\vec 1,\, \} $ obtenemos la siguiente representaci'on matricial:
 \[
\begin{array}{lccr}
&X( 1\, 2\, 3)=
\begin{pmatrix}
1 & \phantom{-}0 &\phantom{-}0 \\ 
0 &\phantom{-}1 &\phantom{-}0 \\
0 &\phantom{-}0 &\phantom{-}1
\end{pmatrix} %%
&X( 2\, 1\, 3 )=
\begin{pmatrix}
1 &\phantom{-}0 &\phantom{-}0 \\ 
0 &-1 &-1 \\
0 &\phantom{-}0 &\phantom{-}1
\end{pmatrix} \\%%% nueva linea
&X(3\, 2\, 1)=
\begin{pmatrix}
1 &\phantom{-}0 &\phantom{-}0 \\ 
0 &\phantom{-}1 &\phantom{-}0 \\
0 &-1 &-1
\end{pmatrix} %%
&X( 1\, 3\, 2 )=
\begin{pmatrix} 
1 &\phantom{-}0 &\phantom{-}0 \\ 
0 &\phantom{-}0 &\phantom{-}1 \\
0 &\phantom{-}1 &\phantom{-}0
\end{pmatrix}\\%%% nueva linea
&X(2\, 3\, 1 )=
\begin{pmatrix} 
1 &\phantom{-}0 &\phantom{-}0 \\ 
0 &-1 &-1 \\
0 &\phantom{-}1 &\phantom{-}0
\end{pmatrix} %%%
&X( 3\, 1\, 2 )=
\begin{pmatrix}
1 &\phantom{-}0 &\phantom{-}0 \\ 
0 &\phantom{-}0 &\phantom{-}1 \\
0 &-1 &-1
\end{pmatrix}
\end{array}
\]
Con respecto a esta segunda base, construida a partir de los  subespacios  invariantes de $V$, las matrices correspondientes  a la representación $X$ son de
la forma
\[
X=
\begin{pmatrix}
A & 0  \\ 
0 & B
\end{pmatrix}
= A \oplus B.
\]
dando una descomposici'on de $X$ como la suma directa de dos representaciones irreducibles: $X= A \oplus B$. En este caso, $A$ es la
representaci'on trivial y $B$ es la representaci'on est'andar. En general, podemos descomponer una representaci'on $X$ de un
grupo finito $G$ como suma de dos representaciones siempre que
exista un subespacio invariante no trivial. Para esto utilizamos el procedimiento que acabamos de ilustrar en la construcci'on de  la representaci'on
est'andar.

\end{Ejercicio}

\begin{definicion}[Representación irreducible] Sea $G$ un grupo finito.
Decimos que una representación de $G$ en $V$ es irreducible si $V$ no tiene ningún subespacio  invariante no trivial bajo la acción de $X$. (Los 
 subespacios  triviales de $V$ son $\{0\}$ y $V$.)
\end{definicion}

La descomposici'on de la representaci'on definici'on de $\S_3$ obtenida en el Ejercicio \ref{descirre} es de hecho la mejor que se puede obtener. Las dos subrepresentaciones que aparecen son irreducibles. En el Ejercicio \ref{standardirred} veremos que la representaci'on est'andar es irreducible. Por otra
parte, al ser unidimensional, es obvio que la representaci'on trivial es irreducible.

%\todo{Lema de Schur}

Si una representaci'on $X$ de $G$ tiene un subespacio no trivial, entonces su complemento ortogonal con respecto a este producto escalar
invariante tambi'en es invariante y no trivial. La restricci'on de $X$ a estos dos subespacios nos produce una descomposici'on m'as
fina. Iterando este proceso, podemos escribir cualquier representaci'on de $G$ como suma directa de representaciones irreducibles.

Un resultado de Maschke nos asegura que cualquier representaci'on compleja de un grupo finito puede escribirse como la suma directa de 
representaciones irreducibles. Esto se demuestra utilizando que para cualquier espacio vectorial complejo  sobre el que  act'ua un grupo finito $G$,
 podemos construir un producto escalar invariante bajo esta acci'on a trav'es del operador de Reynolds: 
 \[
 \langle u,v \rangle = \frac{1}{|G|} \sum_{g \in G} \langle  g \cdot u ,  g \cdot v  \rangle.
 \] 
Por otra parte, si trabajamos con grupos infinitos este resultado no es necesariamente cierto. 
Un grupo que posee la propiedad de que todas sus 
representaciones pueden ser escritas como la suma directa de sus representaciones irreducibles se denomina {\em semi-simple.}

\begin{Ejercicio}
Si $X : G \to GL(V)$ y $Y : G \to GL(W)$ son representaciones de un grupo finito $G$, demostrar que tanto $X\oplus Y$, como $X \otimes Y$ tambi'en lo son.

Demostrar que $X$ induce una representaci'on $X^*$ sobre el espacio dual $V^*$. Si $F \in V^*$ 
\[
X^*(g)  = \phantom{n}^tf(X(g)^{-1}) : V^* \to V^*
\]

Demostrar que los especios $\wedge^k V$ y $Sym^k V$ heredan una estructura de representaci'on de $G$ de aquella de $V$.
\end{Ejercicio}

\begin{Ejercicio}
Si $X$ y $Y$ son representaciones de un grupo finito $G$, demostrar que
\begin{align*}
\chi_{X \oplus Y} &= \chi_X + \chi_Y\\
\chi_{X \otimes Y} &= \chi_X \cdot \chi_Y \\
\chi_{V^*}&= \overline{ \chi_V}\\
\chi_{\wedge^2 V}(g) &= \frac{1}{2} \big[ \chi_V(g)^2-\chi_V(g^2)       \big]
\end{align*}
\end{Ejercicio}

Queremos estudiar la estructura algebraica que poseen las representaciones de un grupo finito. Empezamos por definir un producto interno sobre el espacio de caracteres.
Haremos un breve recuento de algunos hechos de naturaleza general para motivar la definici'on de un producto escalar en el espacio de car'acteres.
Primero observese que del isomorfismo natural entre $Hom(V,W) \cong V^* \otimes W$ deducimos que
\[
Hom_G(V, W) = Hom(V,W)^G \cong (V^* \otimes W)^G
\]
donde $G$ act'ua sobre $\phi \in Hom(V, W)$ como $g(\phi)=g \circ \phi \circ g^{-1}$. (Para entender esta definici'on basta dibujar
 el diagrama conmutativo correspondiente.)

Ahora dada cualquiera proyecci'on $\pi$ de $G$-m'odulo  $U$ sobre su espacio de invariantes $U^G$, la dimensi'on de $U^G$ viene dada por 
\[
tr \, \pi = \frac{1}{|G|} \sum tr(g)
\]
En el caso particular que nos concierne, la proyecci'on de $V^* \otimes W$ sobre $(V^* \otimes W)^G$, tenemos que
\[
tr \, \pi = \frac{1}{|G|} \sum_g \chi_{V^*}(g) \chi_W(g) = \frac{1}{|G|} \sum_g \chi_V(g^-1) \chi_W(g)
\]
ya que si es f'acil ver que su $V$ y $W$ son representaciones, entonces $\chi_{V\otimes W}=\chi_V \chi_W$.  M'as a'un, 
en general $\chi_{V^*}=\overline{\chi_V}$ y por ser $G$ un grupo finito  todos los autovalores son ra'ices de la unidad y $\bar\chi_V= \chi_V(g^{-1})$.

Concluimos que 
\[
 \dim(Hom_G(V_{},W_{})) = \frac{1}{|G|} \sum_{g \in G} \chi_V(g) \chi_W(g^{-1})
\]

\begin{definicion}[Producto interno de caracteres]
Sean $\chi$ y $\phi$ los caracteres correspondientes a dos representaciones del grupo finito $G$. El producto interno entre $\chi$ y $\phi$ se define como:
\begin{align}\label{mult-eq}
\langle \chi, \phi \rangle_G = \frac{1}{|G|} \sum_{g \in G} \chi(g) \phi(g^{-1})
\end{align}
Cuando el grupo $G$ puede ser deducido del contexto, denotamos este producto interno por $\langle \chi, \phi \rangle$.
\end{definicion}

\begin{Ejercicio}
Demostrar que los caracteres de las representaciones trivial y la representaci'on est'andar de $\S_3$ son ortonormales
\end{Ejercicio}

\begin{Teorema}[Relaciones entre caracteres]
Si $X$ y $Y$ son representaciones irreducibles de un grupo finito $G$ con caracteres $\chi$ y $\phi$, tenemos entonces que
\[
\langle \chi, \phi \rangle = \delta_{\chi, \phi}
\]
donde $\delta$ denota la función delta de Kronecker. 
\end{Teorema}

\begin{Ejercicio}[Propiedades de los caracteres] 
Utilizar el teorema anterior para demostrar que si $X$ es una representación matricial de un grupo finito $G$ con car\'acter $\chi$ tal que 
\begin{align}\label{multiplicidades}
X \cong m_1 X^{(1)} \oplus m_2 X^{(2)} \oplus \cdots \oplus   m_k X^{{(k)}},
\end{align}
donde las $X^{(i)}$ son representaciones irreducibles no isomorfas dos a dos, con car\'acter $\chi_i$.
Se tiene entonces que
\begin{enumerate}
\item $\chi = m_1 \chi_1 + m_2 \chi_2 + \ldots + m_k \chi_k.$
\item $\langle \chi, \chi^{(j)} \rangle = m_j   \text{ para cada $j$}.$
\item   $\langle \chi, \chi \rangle =  m_1^2+ m_2^2+\ldots+m_k^2.$
\item $X \text{ es irreducible si y solamente si }  \langle \chi, \chi \rangle =  1.$
\end{enumerate}

\begin{Ejercicio}[La representaci'on est'andar es irreducible] \label{standardirred}
Demuestrar que la representaci'on est'andar es irreducible.

\end{Ejercicio}

Los enteros no-negativos $m_i$ que aparecen en la ecuacion (\ref{mult-eq}) se denominan {\em multiplicidades.} En 
esta situaci'on decimos que la representaci'on
irreducible $X^{(i)}$ aparece con multiplicidad $m_i$ en la representaci'on $X$.

\end{Ejercicio}

%\todo{mover lo de grupo de Grothendieck}

\begin{definicion}[El grupo de Grothendieck de $G$] Sea $G$ un grupo semi-simple.
El {\em grupo de Grothendieck de $G$} es el grupo abeliano libre generado por las clases de isomorf'ia de representaciones irreducibles de $G$ con la operación de suma directa.

Los dos ejemplos can'onicos de grupos semi--simples son los grupos finitos, y el segundo protagonista de nuestras lecciones : el grupo lineal
general.
\end{definicion}

  %Sea $G$ un de ellos. 

%El grupo de Grothendieck de $G$ tiene una estructura algebraica mucho más rica que la de grupo abeliano. Para describirla introducimos las nociones de restricción e inducción de representaciones.

\subsection{\textsf{Restricci'on e inducci'on de representaciones.}}

 Sea $H$ un subgrupo de $G$. Queremos  obtener una representación de $H$ a partir de una representación de $G$ (y viceversa). En la primera situación, el procedimiento es trivial. Si $X$ es una representación matricial de $G$ con car\'acter $\chi$, definimos {\em  la restricción de $X$ a $H$,} que denotamos por $X\downarrow^G_H$, como $X\downarrow^G_H (h)=X(h)$, para cada $h \in H$.  Al car\'acter de la representación resultante lo denotamos por $\chi\downarrow^G_H$. Es importante mencionar que, incluso si $X$ es una representaci'on irreducible $X\downarrow^G_H$, en general,
no lo es.

Construir una representación del grupo $G$ a partir de una representación del subgrupo $H$ es más sutil. Veamos primero que sucede en el caso de un 
grupo finito $G$. Consideremos una transversal 
 $t_1, \ldots, t_k$ de $H$ en $G$; esto es,  una colecci'on de elementos de $G$ tales que
 el grupo $G$ es la unión disjunta de las clases $t_iH$. (Si $G$ es finito $k=|G|/|H|$ tambi'en lo es.)
 
 Sea $Y$ una representación matricial de $H$ sobre el espacio vectorial $V$. 
 Para cada $t_i$, consideramos una copia $t_iV$ de $V$, a cuyos elementos llamamos $t_iv$, con $v \in V$.
 Entonces, el grupo $G$ actúa de manera natural sobre $\oplus_i \, ( t_i Y )$: Para determinar $g(t_iv)$ (donde $g\in G$, $i\in [k]$, $t_iv \in t_iV$) hacemos
 los siguiente: Tenemos que $gt_i \in t_lH$ para un \'unico valor de $l$; sea $gt_i=t_lh$. Entonces, para cada $t_iv\in t_iV$, definimos
 \[
 g \cdot (t_iv) = t_l (h\cdot v) \in t_lV.
 \]
 Se puede demostrar que la representaci'on que resulta no depende de la transversal utilizada para su construcci'on.
 Si calculamos la matriz correspondiente, obtenemos la siguiente definici'on de  {\em la inducción de $H$ a $G$}

 \begin{definicion}[Representación inducida] Sea $H$ un subgrupo de un grupo finito $G$ y sea $Y$ una representación matricial de $H$.
 La matriz correspondiente a la inducción de la representación $Y$ de $H$ a una representación de $G$
 se denota por $Y\uparrow_H^G$ y viene dada por la siguiente matriz por bloques:
 \begin{align*}
  \begin{pmatrix}
 Y(t_1^{-1} g t_1)&Y(t_1^{-1} g t_2)&\ldots &Y(t_1^{-1} g t_k)\\
 Y(t_2^{-1} g t_1)&Y(t_2^{-1} g t_2)&\ldots &Y(t_2^{-1} g t_k)\\
 \ldots&\ldots&\ldots&\ldots\\
 Y(t_k^{-1} g t_1)&Y(t_k^{-1} g t_2)&\ldots &Y(t_k^{-1} g t_k)\\
 \end{pmatrix}
\end{align*}
 donde $Y(g)$ es cero cada vez que $g \notin H$.

  Si $\chi$ es el car\'acter de $Y$, al car\'acter de la representación inducida $Y\uparrow_H^G$ lo denotamos por $\chi \uparrow_H^G$.
 \end{definicion}     
  
\begin{Ejercicio}
 Sea $G=\S_3$ y $H=\{\epsilon, (2,3)\}$,  tenemos que $
G=H \cup (1,2) H \cup (1,3) H$.
Sea $Y=1$ la representaci'on trivial de $H$. Entonces $X=1\uparrow^G_H$ se calcular como sigue: Para construir la primera fila de $X(2\, 1\, 3)$ tenemos que
\begin{align*}
&&Y(\epsilon^{-1}(1, 2) \epsilon)=Y(1,2)=0 &&\text{ya que $(1,2)\notin H$}\\
&&Y(\epsilon^{-1}(1, 2) (1,2))=Y(\epsilon)=1 &&\text{ya que $(1,2)\in H$}\\
&&Y(\epsilon^{-1}(1, 3) \epsilon)=Y(1,3,2)=0 &&\text{ya que $(1,3,2)\notin H$}
\end{align*}
Continuando de esta manera obtenemos que
\[
X( 2\,1\,3)=
\begin{pmatrix}
0&1&0\\
1&0&0\\
0&0&1
\end{pmatrix}
\]
Calcule  $X=1\uparrow^G_H$. Verificar que obtenemos a la representaci\'on definici'on de $\mathbb{S}_3$ (Comparar con el Ejercicio 2.1.) (Tomado de \cite{Sagan})
\end{Ejercicio}
  
  El ejemplo anterior nos ilustra que, a pesar de que la representaci'on trivial es irreducible, $1\uparrow_H^G$ no lo es en general.
  
 Existe una elegante relación entre los procedimientos de inducción y restricción de representaciones.

 \begin{Teorema}[Fórmula de reciprocidad de Frobenius]
 Sea $H$ un subgrupo de un grupo finito $G$,  y sean $\phi$ y $\chi$ car\'acteres de $H$ y de $G$ respectivamente. Tenemos entonces que
 \[
 \langle \phi\uparrow^G_H, \chi \rangle_G = \langle \phi, \chi\downarrow_H^G \rangle_H
 \]
 
 \end{Teorema} 
 
 \subsection{\textsf{$G$--m'odulos.}}
Hasta ahora, siempre hemos trabajados con un espacio vectorial $V$, junto con una de sus bases. Esto nos ha permitido, por ejemplo, asociarle a
cada homomorfismo una matriz. Pero, como  sucede con frecuencia, es m'as f'acil y elegante trabajar con espacios vectoriales
sin la necesidad de fijar una base de antemano. Para esto introducimos  la noción de $G$--módulo.  

Sea $G$ un grupo (no necesariamente finito) y sea $V$ un espacio vectorial. Decimos que $V$ es un {\em $G$-módulo} si existe un homomorfismo de grupos $\rho : G \to GL(V)$ de $G$ al grupo de las transformaciones lineales invertibles del espacio vectorial $V$, $GL(V)$. N'otese que al fijar una base para $V$, cada uno de estos 
homomorfismo $\rho$ nos define   representaci'on matricial de $G$.  Similarmente, cada representaci'on matricial de $G$ nos produce un tal $\rho$.
Las nociones de $G$-m'odulo y de representaci'on matricial de $G$ son equivalentes desde este punto de vista.  Hablaremos de una representaci'on
cuando no queremos enfatizar si estamos trabajando con (o sin) una base.

      Sean $V$ y $W$ dos $G$-módulos. Un {\em $G$-homomorfismo} es una transformación lineal $\theta : V \rightarrow W$ tal que $\theta ( \pi v) = \pi \theta (v)$ para cada $\pi \in G$. Si $\theta$ es invertible, decimos  entonces que es un  $G$-isomorfismo.
 A lo largo de estas lecciones asumimos que todos los $G$-módulos con los que trabajaremos son  de dimensi'on finita.

Denotamos por   $\C[G]$ al 'algebra del grupo $G$, es decir al 'algebra definida por las combinaciones lineales de elementos de $G$, donde la
multiplicaci'on viene dada por  el
producto de $G$. Obsérvese que $\C[G]$ tiene naturalmente  estructura de $H$-m'odulo,  lo que nos permite tomar el siguiente
producto tensorial $\C[H]$-m'odulos. 

%\todo{seguir verificando cuando trabajamos con un grupo finito y cuando no}

La inducci'on de representaciones se puede definir de manera concisa en el lenguaje de $G$-módulos. Sea $H$ un subgrupo de $G$, y sea $V$ un $H$-m'odulo, definimos 
\[
Ind_H^G \, V = \mathbb{C}[G] \otimes_{\mathbb{C}[H]}V
\]
N'otese que $\mathbb{C}[G]$ tiene una estructura natural de $\mathbb{C}[G]$ m'odulo (dada por la restricci'on), y que estamos tomando el producto tensorial
de $\mathbb{C}[H]$--m'odulos.

Hemos obtenido entonces una representación de $G$ a partir de la representaci'on $V$ de $H$ ya que $Ind_H^G \, V$ tiene estructura de $G$--m'odulo. 
Dejamos los detalles de verificar que ambas definiciones coinciden al lector familiarizado con la teor'ia de m'odulos.

\begin{Teorema} Sea $G$ un grupo finito, y sea $\{ V^{(i)} : i \in I\}$ una lista completa de clases de isomorf'ia de $G$--m'odulos. 
Si $\C[G]$ se descompone como $\oplus_{i\in I} m_i V^{(i)},$ tenemos entonces que
\begin{enumerate}
\item[1.] $m_i=\dim V^{(i)}$, para cada $i\in I$.
\item[2.] $\sum_{i\in I} (\dim V^{(i)})^2 = |G|$
\item[3.] El número de clases de conjugación de $G$ es igual al número de representaciones irreducibles de $G$. (Esto es, a la cardinalidad de $I$.)
\end{enumerate}

\end{Teorema}

\begin{Ejercicio} \label{completoS3yS4}
Utilice el  teorema anterior para demostrar que el conjunto de representaciones dado por la representación trivial, la representación alternante, y la  
representación est'andar es un conjunto completo de representaciones irreducibles de $\S_3$.  
\end{Ejercicio}

\subsection{\textsf{Construcci'on de las representaciones del grupo lineal general a partir de las representaciones del grupo sim'etrico.}}

Sea $V$ un espacio vectorial complejo de dimensi'on $m$.	Una representaci'on de $GL(E)$ se dice polinomial si la aplicaci'on
\[
X : GL(V) \to GL(W)
\]
viene descrita por polinomios. Esto es, si después de escoger bases para $V$ y $W$ las $N^2$ funciones coordinadas son polinomios
en las $m^2$ variables determinadas por las entradas de una matriz gen'erica en $GL(V)$. Similarmente, decimos que la representaci'on
es racional u holom'orfica cuando est'an funciones lo son.

En este peque\~no apartado procedemos a dar una construcci'on que nos permite asociar a cualquier representaci'on del
grupo sim'etrico, una representaci'on polinomial del grupo lineal general. Los detalles se encuentran en el libro de William Fulton  \cite{F}.
%A partir de ahora, y al final de cada secci'on, haremos una pausa para discutir las consecuencias de los resultados expuestos en el estudio de la teor'ia de
%representaciones del grupo lineal general. En este momento, veremos como construir 
%representaciones del grupo lineal general, a partir de las representaciones del grupo sim'etrico. Para esto nos basamos en el libro de Fulton \cite{F}.
%% Esta construcci'on nos proporcionar'a un conjunto completo de representaciones irreducibles del grupo lineal general, a partir de la familia de representaciones irreducibles del grupo sim'etrico que estudiaremos en la pr'oxima seccion.
Recordemos que estamos trabajando sobre un espacio vectorial $V$, de dimensi'on $m$. El grupo sim'etrico $\S_n$ act'ua sobre $V^{\otimes n}$ por la derecha
(es decir, permuta las posiciones)
\begin{align*}
&V^{\otimes n} = V \otimes_{\C} V \otimes_{\C} \cdots \otimes_{\C} V,\\
&(u_1 \otimes u_2 \otimes \cdots \otimes u_n) \cdot \sigma = u_{\sigma(1)} \otimes u_{\sigma(2)} \otimes \cdots \otimes u_{\sigma(n)}.
\end{align*}
para cada $u_i \in V$ y $\sigma \in \S_n$. Esta acci'on le proporciona a $V^{\otimes n}$ la estructura e $\S_n$-m'odulo.
Por otra parte, a cada representaci'on $M$ de $\S_n$ le asociamos el espacio vectorial
\[
\mathbb{V}(M) = V^{\otimes n}  \otimes_{\C[\S_n]} M
\] 
En consecuencia tenemos que  $\mathbb{V}(M)$, para cada $w \in V^{\otimes n} $,  $v \in M$ y $\sigma \in \S_n$:
\[
(w \cdot \sigma) \otimes v = w \otimes (\sigma \cdot v)
\]

Como el grupo lineal general $GL(V)$ act'ua a la izquierda en $V$, y como esta acci'on se puede extender diagonalmente  a $V^{\otimes n}$ 
de la siguiente manera 
\[
g \cdot ( u_1 \otimes u_2 \otimes \cdots \otimes u_n ) = g \cdot u_1 \otimes g \cdot u_2 \otimes \cdots \otimes g \cdot u_n.
\]
$V^{\otimes n}$ tambi'en tiene la estructura de $GL(V)$--m'odulo.
Es inmediato ver que ambas acciones conmutan, de manera que $GL(V)$ tambi'en act'ua sobre $\mathbb{V}(M)$ :
$g \cdot (w \otimes v) = (g \cdot w) \otimes v$.

Concluimos entonces que tanto $\mathbb{V}(M)$ como $V^{\otimes n}$  tienen la estructura de $(GL(E),\S_n)$-m'odulo.

Veamos algunos
ejemplos de esta construcci'on. 
Si $M$ es la representaci'on trivial de $\S_n$, entonces $\mathbb{V}(M)$ corresponde a  las potencias sim'etricas de $V$: $\mathbb{V}(M)=Sym^n(V).$ 
Similarmente, si $M$ es la representaci'on alternante
de $\S_n$, obtenemos entonces las potencias sim'etricas de $V$, $\mathbb{V}(M)=\bigwedge^n V$.
Finalmente, si $M=\C[\S_n]$ es la representaci'on regular de $\S_n$, entonces $\mathbb{V}(M)=V^{\otimes n}$.

Esta construcci'on es functorial, dado cualquier homomorfismo $\phi : M \to N$ de $\S_n$--m'odulos, obtenemos un homomorfismo $\mathbb{V}(\phi)= \mathbb{V}(M) \to \mathbb{V}(N)$ de 
$GL(V)$--m'odulos. Una descomposici'on en suma directa $M=\oplus M_i$ de $\S_n$--m'odulos determina una descomposici'on  de $GL(V)$-m'odulos
$\mathbb{V}(M) = \oplus \mathbb{V}(M_i)$.

%%%%%%%%%%%%%%%%%%%%%%%%%%%%%%%%%%%%%%%%%%%%%%%%%
%%%%%%%%%%%%%%%%%%%%%%%%%%%%%%%%%%%%%%%%%%%%%%%%%
%%%%%%%%%%%%%%%%%%%%%%%%%%%%%%%%%%%%%%%%%%%%%%%%%
%%%%%%%%%%%%%%%%%%%%%%%%%%%%%%%%%%%%%%%%%%%%%%%%%figggg

\section{\textsf{Dos familias de representaciones del grupo simétrico.}}

%\todo{Dar una descripcion mas functorial}

El teorema de Maschke nos dice que   cualquier  representación compleja del grupo simétrico puede ser descompuesta como  suma directa
 de representaciones irreducibles. Equivalentemente, cualquier $\S_n$-m'odulo se puede escribir como suma de $\S_n$-m'odulos
irreducibles. Las representaciones irreducibles de $\mathbb{S}_n$ ($\S_n$-m'odulos irreducibles) corresponden a las clases de conjugaci'on del 
 grupo sim'etrico  $\mathbb{S}_n$, que a su vez, como vimos en el ejercicio \ref{index}, vienen indexadas por las particiones de $n$.
 
 En esta secci'on veremos como varios  objetos combinatorios (tableau de Young, tableau semiest'andar, etc) aparecen de
 manera natural en el estudio de estas representaciones.
 Presentamos una  elegante construcción de un conjunto completo de representaciones irreducibles del grupo simétrico 
basada en las notas de Fran\c{c}ois Bergeron \cite{B}.  A partir de este conjunto de representaciones irreducibles complejas del
grupo sim'etrico construiremos una familia completa de representaciones polinomiales complejas irreducibles del grupo lineal general ilustrando
la dualidad que existe entre ambas teor'ias.

Comenzamos por definir una importante familia $\{\H^{\lambda}\}$ de representaciones  por permutaciones de $\S_n$, indexadas  por particiones de $n$. En
general, las representaciones  $H^{\lambda}$ no son irreducibles.  Sin embargo, tienen la siguiente propiedad : Al descomponer a
$\mathcal{H}^\mu$ como suma de representaciones irreducibles, la representaci'on $S^{\mu}$  aparece con multiplicidad uno, y s'olo aquellas  representaciones irreducibles $S^{\lambda}$ indexadas por particiones $\lambda \dom \mu$, en el orden de dominanci'on, aparecen con multiplicidad mayor o igual a uno.
% De Concini y Procesi demostraron que estas representaciones aparecen de manera natural al estudiar las matrices nilpotentes\footnote{  
%Sea $J_{\lambda}$ el ideal asociado a  la clausura de la órbita de matrices nilpotentes con forma de Jordan descrita por $\lambda'$. Sea $I_{\lambda}$ la intersección de $J_{\lambda}$ con el anillo de  las funciones polinomiales sobre las matrices diagonales. De manera equivalente, $I_{\lambda}$ se obtiene al hacer $x_{ij}=0$ en $J_{\lambda}$ cada vez que $i \neq j$
% De Concini y Procesi demostraron que  $R/I_{\lambda}$ tiene la estructura de un $\S_n$--módulo y es isomorfo a la inducción de la representación trivial del grupo de Young $\S_{\lambda}$ a $\S_n$, esto es a $\mathcal{H}^{\lambda}$. 
%},   \cite{DC-P, We}.

Posteriormente, construiremos un conjunto de representaciones irreducibles del grupo simétrico. Concluimos esta sección con la descripción de la  descomposición de las representaciones $\{\H^{\lambda}\}$ como suma directa de representaciones irreducibles.

\begin{definicion}[Tableau  inyectivo] 
Sean $k$  y $n$ dos números dados tales que $k \ge n$.
Un {\em tableau  inyectivo} $t$ de forma $\lambda \vdash n$ es una función inyectiva $t$ del conjunto de casillas del diagrama de  $\lambda$ al conjunto $ [k]$, que denotamos por $t : \lambda \to [k]$.  %Esto es, un tableau en que no utilizamos dos veces el mismo n'umero.
Equivalentemente, es una manera de asociar a cada casilla del diagrama de $\lambda$ un número entre $1, 2, \cdots, k$, sin repetición. 

Sea $l$ una casilla del diagrama de $\lambda$, $t(l)$ la entrada correspondiente a $\lambda$ y $f(l)$ la fila a la que pertenece $l$. Definimos el
peso de un tableau inyectivo $\lambda$ como
\[
x^t=\prod_{l \in \lambda} x_{t(l)}^{f(l)-1}.
\]
Es importante no confundir esta noci'on con la definici'on de peso de un tableau dada en la introducci'on.

 \end{definicion}

\begin{figure}[h]
\begin{eqnarray*}
\begin{matrix} \Young{ 
	\Box{1}& \Box{5}& \Box{4} \cr  
	\Box{7}& \Box{3} \cr 
	\Box{10} \cr}
     \end{matrix} 
 \end{eqnarray*}
\caption{Un tableau inyectivo de forma $\lambda=(3,2,1)$ y peso de tableau inyectivo $x_1^0  x_5^0  x_4^0\,  x_7^1  x_3^1 \, x_{10}^2 = x_3x_7x^2_{10}$. Su estabilizador es 
$\S_{\{1,4,5\}}\times\S_{\{3,7\}}\times\S_{\{10\}}$.}
\label{inyectivo}
\end{figure}

\begin{definicion}[La acción del grupo simétrico sobre el conjunto de los tableaux]
El grupo simétrico $\S_k$ actúa sobre cualquier tableaux  
$t : \lambda \to [k]$ de la manera canónica:
\begin{align*}
\sigma \cdot t  &: \lambda \to [k] \\
\sigma \cdot t(c) &= \sigma(t(c)).
\end{align*}

El estabilizador de esta acci'on sobre $T$ es un {\em subgrupo de Young}, $ \S_{F_1} \times \S_{F_2}  \times \cdots \S_{F_{\ell}},$ donde 
$\S_{F_i}$ es el grupo de permutaciones de las entradas de la $i$-'esima fila de $T$. 

 \end{definicion}

\begin{definicion}
Definimos la familia de $\S_n$-módulos indexados por particiones $\lambda$ de $n$ 
\[
\H^{\lambda} = \L[ x^t : t \text{ es un tableau  inyectivo de forma $\lambda$}\vdash n ]
\]
Nótese que los $\H^{\lambda} $ son representaciones permutación.
\end{definicion}

\smallskip

Por ejemplo, 
$\H^{
(2,1)
}=
\L[ x_3, x_2, x_1]
$
ya que
\begin{align*}
&x^{
\begin{matrix}
\Young{
\Box{1} & \Box{2}  \cr
\Box{3} \cr
}
\end{matrix}}
=
x^{
\begin{matrix}
\Young{
\Box{2} & \Box{1}  \cr
\Box{3} \cr
}
\end{matrix}}
=x_3,\\ %%%%%%%%%%%%%%%%
&x^{
\begin{matrix}
\Young{
\Box{1} & \Box{3}  \cr
\Box{2} \cr
}
\end{matrix}}
=
x^{
\begin{matrix}
\Young{
\Box{3} & \Box{1}  \cr
\Box{2} \cr
}
\end{matrix}}
=x_2,\\%%%%%%%%%%%%%%%
&x^{
\begin{matrix}
\Young{
\Box{3} & \Box{2}  \cr
\Box{1} \cr
}
\end{matrix}}
=
x^{
\begin{matrix}
\Young{
\Box{2} & \Box{3}  \cr
\Box{1} \cr
}
\end{matrix}}
=x_1.
\end{align*}
Como vemos en este ejemplo, en general varios tableaux inyectivos tienen el mismo peso.  Para evitar esta redundancia,
 a\~nadimos  la condición adicional de que las entradas en cada fila se encuentren ordenadas en 
 orden creciente.

\begin{Ejercicio}[Las representaciones $\H^{\lambda}$ son cíclicas]
Demostrar que para cada $\lambda$ fijo  
\[
\H^{\lambda} = \L [\sigma \cdot x^{t(\lambda)} : \sigma \in \S_n].
\]
donde $t(\lambda)$ es el tableau inyectivo can'onico que se obtiene al rellenar las entradas de $\lambda$ con los n\'umeros 
$1, 2, \cdots, \ell(\lambda)$ escribiendo de la manera usual (escribiendo de izquierda a derecha y de arriba hacia abajo.).
\end{Ejercicio}

\begin{Ejercicio}[Ejemplos de módulos $\H^{\lambda}$] Verificar que para $\S_3$,
$\H^{(3)}= \L[1 ]$ es la representación trivial, $\H^{(21)}$ es   la representación definición y $\H^{(111)}$ es la representación regular. 
Más generalmente, demostrar que $\H^{{(n)}}, \H^{(n-1,1)}$ y $\H^{(1^n)}$ son  las  representaciones
 trivial, definición y regular de $\S_n$.

\end{Ejercicio}

\begin{Ejercicio}[Dimensión de los módulos $\H^{\lambda}$]
Demostrar que la dimensión de $\H^{\lambda}$ es
\[
\dim \H^{\lambda} = \frac{n!}{\lambda_1!\, \lambda_2! \, \cdots \lambda_k!}.
\]
\end{Ejercicio}

%\subsubsection*{Las representaciones irreducibles del grupo simétrico} 
La noción de tableau  inyectivo tambi'en nos permite  construir una familia completa de representaciones irreducibles de $\S_n$.

\begin{definicion}[Las representaciones irreducibles del grupo simétrico]
Sea $t$ un tableau  inyectivo, definimos
\begin{align*}
S^{\lambda} =\L[ a_t : t \text{ tableau  inyectivo de forma } \lambda ] 
 \end{align*}
donde $a_t$ es el polinomio que se obtiene al antisim'etrizar a $x^{t(\lambda)}$ con respecto a las columnas de diagrama
de Young de $\lambda$. M'as precisamente, 
sea $\C_t$ el subgrupo de Young definido por el conjunto de todas las permutaciones que estabilizan las entradas en las columnas de $t$. 
Tenemos entonces que 
\[
a_t = \sum_{\sigma \in C_t} \sgn(\sigma) \sigma \cdot x^t
\]
\end{definicion}

\begin{Ejercicio} Verificar  que si $t=\begin{matrix}
\Young{
\Box{1} & \Box{2}  \cr
\Box{3} & \Box{4} \cr
}
\end{matrix}$, entonces su estabilizador es $\mathbb{S}_{\{1,3\}}\times\mathbb{S}_{\{2,4\}}$,   y  
$
 a_t=(x_3-x_1)(x_4-x_2).
$ 
Calcular el polinomio que corresponde al otro tableau est'andar de tipo $(2,2)$.

%Demuestre que si la partici'on correspondiente a $t$ tiene alguna parte repetida, entonces $a_t=0$.
\end{Ejercicio}

\begin{Ejercicio}[Estructura de módulo cíclico]
Demuestre que el grupo simétrico actúa como $\sigma \cdot a_t = a_{\sigma t},$
 y que esto implica que las representaciones $S^{\lambda}$ tienen estructura de módulo cíclico: 
\[
S^{\lambda} =  \L[ \sigma \cdot a_{t(\lambda)} :
 \sigma \in {\S}_n]
\]
Una representaci'on ac'iclica siempre es irreducible.
\end{Ejercicio}

\begin{Teorema}[Teorema Fundamental] Hemos construido un conjunto completo de representaciones irreducibles del grupo simétrico. Más precisamente,
\begin{enumerate}
\item Los ${\S}_n$-módulos $S^{\lambda}$ son irreducibles.
\item $S^{\lambda}$ es isomorfo a $S^{\mu}$ si y s\'olo si $\lambda = \mu$.
\item El conjunto $ \{ S^{\lambda} \}$ es un sistema completo de representaciones irreducibles de ${\S}_n$.
\end{enumerate}

\end{Teorema}
%Para una demostración de este resultado ver \cite{B}.

\begin{Ejercicio}
Demuestre que el conjunto de los tableau est'andar de forma $\lambda$ indexan a una base de $\S^{\lambda}$. Concluya entonces que
\begin{align*}
f^{\lambda} &= \dim S^{\lambda}  \\
n! &= \sum_{\lambda\vdash n} (f^{\lambda})^2
\end{align*}
Esta 'ultima identidad se puede demostrar elegantemente utilizando el famoso algoritmo de Robinson-Schensted-Knuth, ver \cite{B, Sagan, EC2}.
Los n'umeros $\frac{(f^{\lambda})^2}{n!}$ definen  una medida sobre el conjunto de las particiones de $n$, 
la medida de Pancherel del grupo sim'etrico.

\end{Ejercicio}

\begin{Ejercicio}[El determinante de Vandermonde]
El polinomio de Vandermonde $\Delta_n$ se define como el determinante de la matriz  de Vandermonde $(x_i^{n-j})_{i,j=1}^{n}$. 

Demuestre que
\[
\Delta_n = a_{(n-1,n-2,\ldots,1,0)} =  \prod_{1 \le i < j \le n} (x_i - x_j ).
\]
\end{Ejercicio}

Sea $S$ un conjunto con cardinalidad $k$. 
En lo que sigue, denotamos por $\Delta_k(S)$ el determinante de Vandermonde asociado a las variables en $S$.
Veamos como podemos utilizar estas funciones para dar una descripci'on alternativa del conjunto completo de representaciones irreducibles de $\mathbb{S}_n$.

\begin{Ejercicio}[Las representaciones irreducibles de ${\S}_3$.]  
Demostrar que la  siguiente lista nos proporciona un conjunto completo de representaciones irreducibles de ${\S}_3$:  La representación trivial: $S^{3}= {\mathcal L}[1],$
la representación est'andar $S^{2, 1} = {\mathcal L}[ x_3- x_1 \, , x_2- x_1 ]$ y la representación alternante: $S^{1, 1, 1} = {\mathcal L}[\Delta_3]$.
Comparar estos resultados con los obtenidos en el ejercicio \ref{descirre}.

\end{Ejercicio}

\begin{Ejercicio}[Las representaciones irreducibles de ${\S}_4$.]

Demostrar que la siguiente lista nos proporciona un conjunto completo de representaciones irreducibles de ${\S}_4$.
\begin{align*}
&S^{4\phantom{1,1,1}}= {\mathcal L}[1]\\
&S^{3, 1\phantom{1,1}} = {\mathcal L}[ x_4-x_1\, , x_3- x_1 \, , x_2- x_1 \, ,x_4- x_2\, \cdots]\\
&S^{2,2\phantom{1,1}} = {\mathcal L}[(x_1-x_3)(x_2-x_4), (x_1-x_2)(x_3-x_4)]\\
&S^{2,1,1\phantom{1}} = {\mathcal L}[ \Delta_3(x_1,x_2,x_3),  \Delta_3(x_1,x_3,x_4), \Delta_3(x_1,x_2,x_4)]\\
&S^{1,1,1,1}= {\mathcal L}[\Delta_4]
\end{align*}

\end{Ejercicio}

%
%\begin{Ejercicio} \label{caracteres}
%Demostrar que la  tabla de caracteres de las representaciones irreducibles de $\S_3$ viene dada por
%\[
%\begin{array}{c|ccc}
%                 &(1^3)&(2,1)&(3)\\ \hline
%S^{(3)}     &1 & 1& 1\\
%S^{(2,1)}  & 3& 1& 0\\
%S^{(1^3)} & 6& 0& 0
%\end{array}
%\]
%Calcular la tabla de caracteres correspondiente a las representaciones irreducibles de $\S_4.$ 
%\end{Ejercicio}

En general, se conoce muy poco acerca de cómo descomponer una representación como suma directa de sus  representaciones irreducibles. El caso particular  de las representaciones de $\H^{\lambda}$ es   especial. Su descomposici'on se describe simplemente  en
el lenguage de la combinatoria de los tableros de Young y viene descrita por la  {\em regla de Young}.
\begin{Teorema}[La regla de Young] \label{Young}
Descomposición de la representación $\H^{\mu}$ en términos de sus componentes irreducibles.
\[
\H^{\mu} = \oplus_{\lambda} K_{\lambda, \mu} S^{\lambda}
\]
donde los $K_{\lambda, \mu}$ son  los  coeficientes de Kostka.
\end{Teorema}

%\todo{Regla de Young ? Regla de Pieri ?}

\begin{Ejercicio}[Un ejemplo de la regla de Young] Demostrar que 
\[
\H^{( 3, 2, 1 )} = S^{(3, 2, 1)} \oplus S^{(3, 3)} \oplus 2 \cdot  S^{(4, 2 )} \oplus S^{(4, 1, 1)} 
\oplus  2 \cdot   S^{(5, 1)} \oplus S^{(6)}.
\]

\end{Ejercicio}

Sabemos que el coeficiente de Kostka $K_{\lambda,1^n}= f^{\lambda}$, es igual al número de  tableaux estándar de forma $\lambda$. Tenemos 
entonces que
\[
\H^{(1^n)} = \oplus_{\lambda \vdash n} ( S^{{\lambda}})^{\oplus f^{\lambda}}
\]

 %%%%%%%%%%%%%%%%%%%%%%%%%%%%%%%%%%%%%%%%%%%%%%%%%
%%%%%%%%%%%%%%%%%%%%%%%%%%%%%%%%%%%%%%%%%%%%%%%%%
%%%%%%%%%%%%%%%%%%%%%%%%%%%%%%%%%%%%%%%%%%%%%%%%%
%%%%%%%%%%%%%%%%%%%%%%%%%%%%%%%%%%%%%%%%%%%%%%%%%

\subsection{\textsf{Las representaciones irreducibles del grupo lineal general.}} En la secci'on anterior vimos como obtener las potencias sim'etricas
y alternantes a partir de las representaciones trivial y alternante del grupo sim'etrico. Veamos qué sucede al aplicar el functor $\mathbb{V}(\cdot)$ a las
representaciones del grupo sim'etrico que acabamos de construir.

Como un ejercicio, el lector deber'a verificar que  dada cualquier partici'on $\lambda=(\lambda_1, \lambda_2, \cdots, \lambda_l)$,
obtenemos la siguiente representaci'on de  $GL(V)$:
\[
\mathbb{V}(\mathcal{H}^\lambda)\cong Sym^{\lambda_1}(V)\otimes Sym^{\lambda_2}(V)\otimes \cdots \otimes Sym^{\lambda_l}(V).
\]
Al igual que en el caso del grupo sim'etrico, esta representaci'on no es irreducible, y su descomposici'on en irreducibles viene dada
por los coeficientes de Kostka.

A partir de las representaciones irreducibles de los grupos sim'etricos $\S_n$, obtenemos a un conjunto completo de representaciones
irreducibles de $GL(V)$.

\begin{Teorema} La familia de representaciones polinomiales del grupo lineal general $GL(V)$
\[
V^{\lambda} := \mathbb{V}(S^{\lambda})
\]
 donde $\lambda$ es una partici'on de longitud $\le \dim(V)$, es una familia completa de representaciones polinomiales irreducibles
 de $GL(V)$. N'otese que, al no tener ninguna condici'on sobre $|\lambda|$ esta es una familia infinita.

M'as a'un, tenemos la siguiente descomposici'on de $V^{\otimes n}$ como $GL(V)$-m'odulo : 
\[ 
V^{\otimes n} = \mathbb{V}(\C[\S_n]) 
\cong \bigoplus_{\lambda \vdash n}
 (V^{\lambda})^{\oplus f^{\lambda}}
\]	
donde la suma se toma sobre todas las particiones de $n$ de longitud menor o igual a $dim(V)$.
\end{Teorema}

El algoritmo de Robinson-Schensted-Knuth (RSK), ver \cite{EC2}, nos proporciona una biyecci'on entre el n'umero de funciones
$[n] \to [m]$, y las parejas de tableaux, el primero estándard con entradas en $n$, el segundo semi-est'andar con entradas en $m$.
En este contexto nos proporciona una demostraci'on combinatoria de esta identidad, al nivel de espacios vectoriales.

%%%%%%%%%%%%%%%%%%%%%%%%%%%%%%%%%%%%%%%%%%%%%%%%%
%%%%%%%%%%%%%%%%%%%%%%%%%%%%%%%%%%%%%%%%%%%%%%%%%
%%%%%%%%%%%%%%%%%%%%%%%%%%%%%%%%%%%%%%%%%%%%%%%%%
%%%%%%%%%%%%%%%%%%%%%%%%%%%%%%%%%%%%%%%%%%%%%%%%%

\section{\textsf{Las funciones simétricas.}}

Ya en nuestra primera lecci'on, introdujimos las funciones sim'etricas a trav'es del estudio de las matrices no negativas.
En esta segunda lecci'on las estudiaremos desde un punto de vista algebraico. Para esto haremos
un 'enfasis particular en la relación que existe entre la estructura de álgebra de Hopf de funciones
 simétricas y la teoría de  representaciones  del  grupo simétrico y del grupo lineal general.

%%%%%%%%%%%%%%%%%%%%%%%%%%%%%%%%%%%%%%%%%%%%%%%%%
%%%%%%%%%%%%%%%%%%%%%%%%%%%%%%%%%%%%%%%%%%%%%%%%%
%%%%%%%%%%%%%%%%%%%%%%%%%%%%%%%%%%%%%%%%%%%%%%%%%
%%%%%%%%%%%%%%%%%%%%%%%%%%%%%%%%%%%%%%%%%%%%%%%%%

%\subsection{El álgebra de las funciones simétricas.} 
A un conjunto ordenado de variables lo llamamos alfabeto. Resulta muy útil escribir un alfabeto como una suma formal de variables. Por ejemplo, 
 el alfabeto $X= \{x_1, x_2, \ldots \}$  lo 
escribimos como 
$
X=x_1+x_2+\cdots
$.
Denotamos por $\mathbb{Z}[|X|]$ al 'algebra de las series formales en el alfabeto $X$.
Si $\pi$ es una permutaci'on en $\mathbb{S}_n$ y $f \in \mathbb{Z}[|X|]$  definimos
\[
  f[X] \cdot\pi =p(x_1, x_2, \cdots, x_n, x_{n+1} \cdots) = 
f(x_{\pi(1)}, x_{\pi(2)}, \cdots, x_{\pi{(n)}},  x_{n+1} \cdots)
\]

Denotamos por $\Sym$ la  subálgebra de $\Z[|X|]$ formada por las series, con un
n'umero finito de componentes no nulas, que son invariantes bajo esta acción (para cada $n \in \mathbb{N}$). Por ejemplo, 
\begin{align*}
p_k=x^k_1+x^k_2+x^k_3+\ldots &\in \Sym\\
   1+p_k+p_k^2+\ldots = \frac{1}{1-p_k} &\notin \Sym
\end{align*}
El segundo es un ejemplo de una funci'on que, a pesar de ser invariante, tiene un n'umero infinito de componentes no nulas. Tales funciones
pertenecen a la completaci'on de $\Sym$.
Denotamos por $\Sym^{(k)}$a la componente homogénea de grado $k$ de $\Sym$.

 \begin{definicion}[Las funciones simétricas monomiales, $m_{\lambda}$]

Dada $\lambda=(\lambda_1, \lambda_2, \cdots, \lambda_{\ell})$, definimos $m_{\lambda}$ como la  suma de todos los  monomios {\em diferentes} que se obtienen al simetrizar  $x^{\lambda_1}_1 x^{\lambda_2}_2 \cdots x^{\lambda_{\ell}}_{\ell}$. Por ejemplo, $m_{2,1,1}$ es igual a
\[
m_{2,1,1}=x_1^2x_2x_3+x_2^2x_1x_3+x_3^2x_1x_2+\cdots
\]
\end{definicion}

\begin{Ejercicio}[El espacio vectorial $\Sym^{(k)}$] \label{basemonomios}
Demuestre que el conjunto de las funciones simétricas monomiales $m_{\lambda},$ con $ \lambda \vdash k$,  constituye una base   para $\Sym^{(k)}$.
En particular,   la dimensión de $\Sym^{(k)}$ viene dada por el número de particiones de $k$. 

\end{Ejercicio}

Sea $\lambda=(\lambda_1, \lambda_2, \cdots, \lambda_{\ell})$ una partición. Definimos tres familias de funciones simétricas, las 
llamadas {\em  bases multiplicativas,}
\begin{align*}
&&e_{\lambda}=e_{\lambda_1} e_{\lambda_2} \cdots e_{\lambda_{\ell}},
&&p_{\lambda}=p_{\lambda_1} p_{\lambda_2} \cdots p_{\lambda_{\ell}}, 
&&h_{\lambda}=h_{\lambda_1} h_{\lambda_2} \cdots h_{\lambda_{\ell}}. 
\end{align*}
donde $h_n$ es igual a la suma de todos los monomios de grado $n$, $e_m$ la suma de todos los monomios de grado $n$ libres de cuadrados,
y $p_n$ la suma de las potencias $n$-'esimas de los elementos de $X$. 
Por ejemplo, si $X=x+y+z$, entonces $h_2[X]=x^2+y^2+z^2+xy+xz+yz$, $e_2[X]=xy+xz+yz$ y $p_2[X]=x^2+y^2+z^2$.

\begin{Teorema}[Teorema Fundamental de las funciones simétricas]\label{fundamental} Sea $X$ un alfabeto. 
La familia de funciones simétricas  $e_{\lambda}[X]$  es una base (como $\mathbb{Z}$--m'odulo) para $\Sym$. 
\[
Sym[X]={\Z}[e_1[X], e_2[X], \cdots].
\] 
 \end{Teorema}

Este teorema no es dif'icil de demostrar. Se puede hacer estudiando la matriz de cambio de bases entre las funciones
sim'etricas elementales y las monomiales, y observando que son triangulares y con todas sus entradas en la diagonal principal
iguales a uno. O se puede
hacer utilizando al orden lexicogr'afico en el algoritmo de la divisi'on de Buchberger.

%Análogamente, los conjuntos $\{p_{\lambda}  \}$ y $\{h_{\lambda}   \}$
%son bases para $\Sym$. (En el caso de las sumas de potencias es necesario tomar los coeficientes sobre $\mathbb{Q}$, ver Ejercicio \ref{powerq})

 En el siguiente ejercicio estudiaremos las matrices de cambio de base entre las bases multiplicativas de $\Sym$.  Obtendremos
 entonces que las tres familias multiplicativas que hemos introducidos son, como su nombre lo sugiere, bases para $Sym$.

\begin{Ejercicio}[Familias multiplicativas, matrices de cambio de base] 
\label{fgmultiplicativas}
Definimos las siguientes funciones generatrices asociadas a las bases multiplicativas de $\Sym$:
\begin{align*}
E(t)&= \sum_{k \ge 0} e_kt^k = \prod_{i\ge1} (1+tx_i)\\
H(t)&=\sum_{k \ge 0} h_kt^k = \prod_{i\ge1}(1-tx_i)^{-1}\\
P(t)&=\sum_{k\ge 0} p_{k+1}t^k
\end{align*}
Demuestre las siguientes identidades entre las funciones generatrices $E(t), H(t)$ y $P(t)$:
\begin{align}\label{multiplicativasfg}
&&E(t)H(t)=1,
&&P(t)=\frac{H'(t)}{H(t)}, &&P(-t)=\frac{E'(t)}{E(t)}.
\end{align}
De estas identidades, deduzca que para cada entero $k>0$
\begin{enumerate} 
\item Recurrencia entre las elementales y las homog'eneas : $\sum_{i+j=k}(-1)^ie_ih_j=0$.

\item Concluya que las funciones sim'etricas homog'eneas constituyen una base para $Sym$.

\item Identidades de Newton: $k p_k = \sum_{i+j=k}p_ih_j=\sum_{i+j=k}(-1)^{j-1}p_ie_j. $ 

\item La funci'on completa homog'enea $h_n$ y las series de potencias : $h_n=\sum_{\lambda \vdash n} z_{\lambda}^{-1} p_{\lambda}$

\item La funci'on elemental $e_n$ y las series de potencias : $e_n=\sum_{\lambda\vdash n}(-1)^{n+\ell(\lambda)} z^{-1}_{\lambda} p_{\lambda}$.

\end{enumerate}

%\todo{Esto estaba inicialmente escrito para $Sym_n$, revisar que lo que quedo sigue teniendo sentido.}

\end{Ejercicio}

%\begin{Ejercicio}[Las sumas de potencias son una base de $\Sym \otimes \Q$] \label{powerq}
% Utilizar la regla de Cramer para obtener la siguiente identidad
%\[
%n!\,e_n = 
%\begin{vmatrix}
%&p_1  &1 &0 &\cdots &0 &0\\
%&p_2 &p_1 &2 &\cdots &0 &0\\
%&\cdots &\cdots &\cdots &\cdots &\cdots &\cdots\\
%&p_{n-1} &p_{n-2} &p_{n-3} &\cdots &p_1 & n-1\\
%&p_n&p_{n-1} &p_{n-2} &\cdots &p_2 &p_1
%\end{vmatrix}
%\]
%Deducir de esta identidad  que las series de potencias $p_{\lambda}$ son una base para el álgebra $\Sym \otimes \Q,$ las funciones
%sim'etricas con coeficientes racionales.
%\end{Ejercicio}
 Esto nos permite definir $\Sym$ como el 'algebra de polinomios $\mathbb{Z}[e_1, e_2, \cdots]$, donde las $e_i$ son variables formales,
 y luego introducir las otras bases de $\Sym$ utilizando las identidades que acabamos de encontrar.

\begin{Ejercicio}[La involuci'on $\omega$]
Definimos un endomorfismo
\begin{align*}
\omega : Sym &\to Sym\\
\omega(e_n) &= h_n
\end{align*}
Demuestre que $\omega$ es una involuci'on utilizando la recursi'on entre las homogeneas y las elementales obtenida en el 
Ejercicio \ref{fgmultiplicativas}.
\end{Ejercicio}

Ya en la primera de nuestra serie de lecciones introdujimos la más intrigante de las bases para el álgebra de las funciones simétricas; la base de Schur. En aquel momento las funciones de Schur se definieron utilizando  la matriz de Jacobi--Trudi que las expresa en la base de las funciones completas homogéneas. Daremos ahora la definición combinatoria para esta importante base.

%El peso de un tableau semi--estándar $T$ de forma $\lambda$ y contenido $\alpha=(\alpha_1, \alpha_2, \ldots)$ se define como
%\[
%x^{T}=x_1^{\alpha_1(T)}x_2^{\alpha_2(T)}\ldots
%\]
%Por ejemplo,

%\begin{figure}[h]
%\begin{eqnarray*}
%\begin{matrix} \Young{ 
%	\Box{1}& \Box{1}& \Box{2} & \Box{2}& \Box{6}\cr  
%	\Box{2}& \Box{2}& \Box{3} &\Box{4} \cr 
%	\Box{4} &\Box{6}\cr}
%     \end{matrix} 
% \end{eqnarray*}
%\caption{Un tableau semi-estándar de forma $\lambda=(5,4,2)$ y contenido $\alpha=(2,4,1,2,0,2)$, tiene peso $x_1^2 x_2^4 x_3 x_4^2 x_6^2$.}
%\label{semistandar}
%\end{figure}
\begin{definicion}[Definición combinatoria de las funciones de Schur.] Sea $X=x_1+x_2+\ldots$ un alfabeto y sea $\lambda$ una partición de $n$. 

La función de Schur $s_{\lambda}$ se define como
\[
s_{\lambda}[X]=\sum_{T}x^T
\]
 donde sumamos sobre todos los tableaux semi--estándar $T$ de forma $\lambda$ y donde $x^T$ denota el peso del tableau $T$.
 
Definimos las funciones de Schur sesgadas $s_{\lambda/\mu}[X]$ similarmente, pero ahora sumamos sobre todos
 los tableaux semi-est'andar sesgados $T$ de forma  $\lambda / \mu$.
\end{definicion}

Por ejemplo, 
\begin{align*}
s_{(2,1)}= m_{(2,1)}+2 m_{(1,1,1)}.\\
s_{(2,1)/(1)}=h_1^2.
\end{align*}
No es evidente de esta definici'on que las funciones de Schur sean sim'etricas. Esto se puede ver
directamente utilizando un elegante argumento de D. Knuth (ver \cite{Sagan}).

\begin{Ejercicio} Utilice la definici'on combinatoria de una funci'on de Schur para demostrar que
\begin{align*}
s_{(n)}&=h_n\\
s_{(1^n)}&=e_n\\
h_{(n-1,1)}&=s_{(n)}+s_{(n-1,1)}.
\end{align*} 
\end{Ejercicio}

De la definición combinatoria de las funciones de Schur se puede conseguir su desarrollo en la base de las funciones simétricas monomiales. Recuerde que el coeficiente de Kostka $K_{\lambda, \mu}$ se define como el número de tableaux semi--estándar de forma $\lambda$ y contenido $\mu$.  
El hecho de que las funciones de Schur sean sim'etricas es equivalente a la igualdad $K_{\lambda,\mu}=K_{\lambda,\tilde\mu}$ para cualquier rearreglo de las parte de $\mu$. 
Concluimos que
\[
s_{\lambda} = \sum_{\mu} K_{\lambda, \mu} m_{\mu}.
\]

\begin{Ejercicio}[Las funciones de Schur son una base para $\Sym$]
Demuestre que el conjunto de las funciones de Schur, $\left\{s_{\lambda}\right\}_\lambda$, donde $\lambda\vdash k$, forma una base para 
$\Sym^{(k)}$.
Para esto demuestre  que el determinante de la matriz de Kostka (la matriz definida por los coeficientes de Kostka en el orden
lexicogr'afico) tiene determinante diferente de cero.
\end{Ejercicio}

\begin{Teorema}[La involuci'on $\omega$ y la base de Schur]
La involuci'on $\omega$, aplicada a la base de Schur, tiene el efecto de transponer la partici'on que le sirve de 'indice
\[
\omega(s_{\lambda})=s_{\lambda'}
\]
\end{Teorema}

\begin{obs}
Con frecuencia, no es demasiado importante el tama\~no del alfabeto con el que trabajamos. (No siempre. N'otese que $s_{1,1,1}[X]$ es cero si $|X|<3.$) 
Por otra parte, para cada una de las bases que hemos definido
\[
b(x_1, \cdots, x_l,0,\cdots,0)=b(x_1,\cdots, x_l)
\]
Esto implica, por ejemplo, las identidades entre las funciones de Schur en $m$ variables, son ciertas para las funciones de Schur, en un
n'umero infinito de variables, cuando imponemos la condici'on adicional que todas las particiones que aparezca tengan longitud $\le m$.
\end{obs}

Definimos un producto escalar sobre $\Sym$, que denotamos por $\langle \,\,, \, \rangle$, diciendo que la base de Schur es una base ortonormal:
\[
\langle s_{\mu}, s_{\nu} \rangle = \delta_{\mu,\nu}.
\]
donde $\delta_{\mu,\nu}$ es la funci'on delta de Kronecker.
De la ortonormalidad de la base de Schur vemos  que la involuci'on $\omega$ es una isometr'ia.

Este mismo producto escalar se puede definir de manera equivalente diciendo que $\langle m_{\lambda}, h_{\mu} \rangle=\delta_{\lambda,\mu},$ o que $\langle p_{\lambda}, p_{\mu} \rangle=\delta_{\lambda,\mu}z_{\lambda}.$  

\subsection{\textsf{La estructura de 'algebra de $Sym$.}}

El producto de dos funciones sim'etricas es una funci'on sim'etrica, de manera que $\Sym$ tiene una estructura de 'algebra. 

%\todo{yadah}

\begin{definicion}[Los coeficientes de Littlewood--Richardson y el producto de funciones de Schur]
Dado que las funciones de Schur son una base para $Sym$, para cualquier $\mu$ y $\nu$ podemos encontrar constantes $c^{\lambda}_{\mu,\nu}$ tales
que
 \[
 s_{\mu}s_{\nu} = \sum_{\lambda} c^{\lambda}_{\mu,\nu} s_{\lambda}
 \]
 Esta identidad se traduce en 
 \[
c^{\lambda}_{\mu,\nu}=\langle s_{\lambda}, s_{\mu}s_{\nu} \rangle
\]
Los coeficientes $c^{\lambda}_{\mu,\nu}$ se conocen como {
\em  los coeficientes de Littlewood-Richardson.}
\end{definicion}

Los coeficientes de Littlewood-Richardson juegan un rol fundamental en la teor'ia de
las funciones sim'etricas (y sus aplicaciones a la teor'ia de representaciones, entre muchas otras 'areas). La famosa regla de Littlewood--Richardson nos
proporciona una interpretaci'on combinatoria para estos coeficientes (cuentan la cardinalidad de un cierto conjunto de tableaux semi-est'andar definido
en funci'on de $\lambda, \mu$ y $\nu$.) Esta interpretaci'on combinatoria tambi'en nos proporciona una algoritmo (ineficiente) para calcularlos.

%\todo{yadda yadda yadda}

Veamos ahora otra instancia en la que los coeficientes de  Littlewood-Richardson aparecen de manera natural en la teor'ia de las funciones sim'etricas.

%\begin{definicion} Sean $\lambda$ y $\mu$ particiones tales que, si identificamos a las particiones con sus diagramas de Young,  $\mu \subseteq \lambda$.
%Definimos a la funci'on de Schur sesgada 
%\[
%s_{\lambda / \mu} =\sum_{T}x^T
%\]
%donde $x^T$ denota al peso del tableau $T$, y donde sumamos sobre todos los tableaux semi--estándar $T$ de forma $ \lambda / \mu$. 

%
%\end{definicion}

\begin{definicion}[Los coeficientes de Littlewood-Richardson y el adjunto de la multiplicaci'on por una funci'on de Schur]
Las funciones de Schur sesgadas se expresan en la base de Schur utilizando los coeficientes de Littlewood-Richardon.
\[
s_{\lambda \setminus \mu} =\sum_{\nu} c^{\lambda}_{\mu,\nu} s_{\nu}
\]
Por lo tanto, las funciones de Schur sesgadas son el adjunto de la multiplicaci'on por una funci'on de Schur.
M'as precisamente, definimos una operaci'on lineal $s_{\mu}^{\perp}$ sobre $Sym$, declarando 
$s_{\mu}^{\perp} s_{\lambda}=s_{\lambda/  \mu}$.  (Observe que por definici'on $s_{\mu}^{\perp} s_{\lambda}=0$ si $\mu$  no se encuentra
contenido en $ \lambda$.)

La identidad 
$
c_{\mu,\nu}^{\lambda} = \langle s_{\lambda}, s_{\mu}s_{\nu} \rangle  = \langle  s_{\lambda / \mu}, s_{\nu}      \rangle
%\langle  s_{\mu}^{\perp}s_{\lambda}, s_{\mu}      \rangle
$ implica que $s_{\mu}^{\perp}$ es el adjunto de la multiplicaci'on por $s_\mu$; es decir, que para cualquier par de funciones sim'etricas $f$ y $g$ :
\[
 \langle s_{\mu}  \,  f,  g \rangle  =  \langle f,  s_{\mu}^{\perp}  g \rangle.
\]
\end{definicion}

\begin{Ejercicio} Calcular el adjunto de la multiplicaci'on por $p_{\lambda}$, $e_{\lambda}$ y $h_{\lambda}$.
\end{Ejercicio}

\section{\textsf{El \'algebra de Hopf de las funciones sim'etricas}}
	
	Las funciones simétricas, junto con la operación de multiplicación, forman un álgebra graduada, pero tienen una estructura mucho más rica. Antes de empezar a describirla introduciremos las nociones de  co'algebra, bi'algebra y 'algebra de Hopf.

\begin{definicion}[Coálgebra]
Sea $\kk$ un anillo conmutativo.
Una {\em coálgebra $C$} es un $\kk$--módulo junto con una {\em comultiplicación $\Delta : C \to C \otimes C$}
y una {\em counidad} $\epsilon: C \to \kk$ que satisfacen las propiedades coasociativa y counitaria. Esto es,
\begin{align*}
(\Delta \otimes 1) \Delta &= (1 \otimes \Delta) \Delta\\
(\epsilon \otimes 1)\Delta &= (1 \otimes \epsilon) \Delta
\end{align*}
 donde $1$ es la identidad de $C$.

\end{definicion}

\begin{definicion}[Biálgebra]
Sea $B$ es un $\kk$--módulo que es simultáneamente un álgebra y una coálgebra. Decimos que es una {\em biálgebra} si ambas estructuras son compatibles. Esto es, si tanto la comultiplicación como la counidad son morfismos (unitarios) de álgebras. Equivalentemente, se puede pedir que tanto la multiplicación  como la unidad sean morfismos de coálgebras.

\end{definicion}

\begin{definicion}[Álgebra de Hopf]

Una aplicación $\kk$--lineal $\psi : H \to H$ sobre una biálgebra $H$ se denomina ant'ipoda
si para cada $h$ en $H$ se tiene que 
\[
\sum h_i \psi(h_i')=\epsilon(h)1= \sum \psi(h_i)h'_i
\]
 donde $\Delta h= \sum h_i \otimes h'_i$.
  Un álgebra de Hopf es una biálgebra que posee  una ant'ipoda.
\end{definicion}

%
%\subsubsection*{La operación de sustitución. Operaciones sobre alfabetos.}

%Alain Lascoux desarrolla la teoría de las funciones simétricas dándole un papel central a las operaciones sobre alfabetos y la operación de sustitución, \cite{Lascoux} . Recomendamos la  lectura de donde, entre otras muchas cosas, se dan aplicaciones de las funciones simétricas a operaciones sobre polinomios como por ejemplo el algoritmo de Euclides, las resultantes y el método de interpolación de Lagrange.

Recuerdese que denotamos por $f[X]$ a la función  $f$ evaluada en el alfabeto $X$. Explícitamente,
\[
f[X]= f(x_1, x_2, \ldots)
\]
Ahora queremos pasar a  considerar una funci'on sim'etrica   como un alfabeto. Por ejemplo, a la funci'on sim'etrica
$
p_2 = x_1^2 + x_2^2 + x_3^2+ \cdots
$ la identicamos con el conjunto de ``variables''  $\{ x_1^2, x_2^2, x_3^2, \cdots \}$. 

\begin{definicion}[El pletismo  de funciones sim'etricas]
Definimos el {\em pletismo (o sustituci'on)}  de funciones sim'etricas trabajando con la base de las series de potencia. Como
nos sugiere el ejemplo anterior, queremos que
$p_n[p_2]=p_{2n}$. Para esto, procedemos de la siguiente manera:
Sea $g$ una serie de potencias formal sobre el anillo de las series formales $\Q[[x_1, x_2, \cdots]]$. Escribimos 
\[
g=\sum_{\alpha}   u_{\alpha}
\]
donde cada $u_{\alpha}$ es un monomio (con coeficiente uno) en $\Q[[x_1, x_2, \cdots]]$. Por ejemplo, $2x=x+x$.

Definimos
\begin{align*}
p_n[g] &=  \sum_{\alpha} c_{\alpha} u_{\alpha}^n\\
p_{\lambda}[g]&= p_{\lambda_1}[g] p_{\lambda_2}[g] \cdots  p_{\lambda_n}[g]
\end{align*}
donde $\lambda=(\lambda_1, \lambda_2, \cdots, \lambda_n)$.  

Por ejemplo, utilizando la definici'on de pletismo, vemos que
\[
p_n[2x+2y]=p_n[x+x+y+y]=x^n+x^n+y^n+y^n=2x^n+2y^n.
\]

Finalmente, definimos $f[g]$, para cualquier función simétrica $f$, diciendo que $f[g]$ es lineal en $f$. Para calcularlo
empezamos por expresar a $f$ en la base de las series de potencias.
\end{definicion}

El pletismo {\em no} es una transformaci'on lineal, por ejemplo $p_{\lambda}[2 x] =2^{\ell(\lambda)}p_{\lambda}[x].$

%\todo{alfabetos negativos y la involuci'on omega}

\begin{Ejercicio}
Demuestre las siguientes propiedades de la operación de sustitución.
\begin{enumerate}
\item
Si $f$ y $g$ son funciones simétricas, entonces $f[g]$ también lo es.
\item La operación de sustitución es asociativa.
\item Se tiene que
 $p_n[p_m] = p_{nm}.$
\item
En general, $f[g] \ne g[f]$.
\item
Para cualquier función simétrica $f$ se tiene que $p_n[f]=f[p_n]$.
\end{enumerate}
\end{Ejercicio}

%	\begin{Ejercicio}

%	Calcular $ h_n[p_m]$, $h_m[e_2]$ y $h_m[h_2]$.
%	\end{Ejercicio}
%	\todo{como está es demasiado difícil el ejercicio}
\begin{definicion}
Dados dos alfabetos $X= x_1+x_2+\cdots$ e $Y=y_1+y_2+\cdots$, definimos los alfabetos suma y producto como
\begin{align*}
X+Y&=  x_1+y_1+ x_2 + y_2 + \cdots \\
XY\phantom{+}&= x_1y_1+x_1y_2+\cdots+x_iy_j+\cdots
\end{align*}
\end{definicion}

\begin{Ejemplo}[Los coeficientes de Littlewood--Richardson y la suma de alfabetos]
Para evaluar $s_{\lambda}[X+Y]$ podemos asumir que todas las variables  $x_i \in X$ son menores que todas las $y_j \in Y$, esto es
posible ya que  $s_{\lambda}$ es una funci'on sim'etrica. 
De la definici'on combinatoria de una  funci'on de Schur sesgada obtenemos que $s_{\lambda}[X+Y]$ es una suma de funciones sim'etricas
de la forma $s_\mu[X]s_{\lambda / \mu}[Y]$, donde
$\mu \in \lambda$ es la partici'on que definen las entradas en $X$ y $\lambda / \mu$ es la partici'on que definen las
entradas en $Y$. 

\begin{align*}
s_{\lambda}[X+Y]&=\sum_{\mu\in \lambda} s_{\mu}[X] s_{\lambda/ \mu}[Y]=\sum_{\mu, \nu} c_{\mu, \nu}^{\lambda} \, s_{\mu}[X] s_{\nu}[Y]
\end{align*}
\end{Ejemplo}

Un problema abierto importante en la combinatoria algebraica es el de entender el desarrollo de la funci'on  $s_{\lambda}[XY]$ 
en la base de Schur. Por otra parte, el desarrollo de $s_{(n)}[XY]=h_n[XY]$ es particularmente 'util y viene dada por el
kernel de Cauchy, que nos proporciona el desarrollo de $h_n[XY]$  en cualquier par de bases duales\footnote{Las bases $u_{\lambda}$ y $v_{\lambda}$ son {\em duales} si 
$\langle u_{\lambda}, v_{\nu}\rangle=\delta_{\lambda,\nu}$}   $u_{\lambda}$ y $v_{\lambda}$ :

\begin{Teorema}[El kernel de Cauchy]
\begin{align*}
h_n[XY] &= \prod_{i,j} \frac{1}{1-x_iy_j}\\
&=\sum_{\lambda} u_{\lambda}[X]v_{\lambda}[Y] \\
&=\sum_{\lambda} s_{\lambda}[X]s_{\lambda}[Y]  
\end{align*}
Adem'as de las funciones de Schur, los otros dos pares mas importantes de bases de funciones sim'etricas ortonormales son el par
$\{h_{\lambda}\}$ y $\{m_{\lambda}\}$,  y el par $\{p_{\lambda}\}$ y $\{p_{\lambda}/z_{\lambda}\}$.
\end{Teorema}

\subsection{\textsf{La estructura de 'algebra de Hopf de $Sym$, la suma de alfabetos y los coeficientes de Littlewood--Richardson.}}

Pasamos ahora a describir la estructura de álgebra de Hopf de $Sym$ con respecto a la multiplicación ordinaria de series de potencias.
Identificamos $\Sym \otimes \Sym$ con las funciones en dos alfabetos $X$ y $Y$ que son simétricas en cada alfabeto separadamente. 
Bajo esta identificaci'on $f\otimes g$ corresponde al producto $f[X]g[Y]$. 

Utilizando la operación de suma de alfabetos, definimos  una operación de comultiplicación sobre $\Sym$,
\begin{align*}
\Delta : \Sym &\to \Sym \otimes \Sym\\
\Delta f &= f[X+Y]
\end{align*}
Se tiene que $(\Sym, \Delta)$ tiene estructura de coálgebra. La counidad $\epsilon$ viene dada por la proyección $ \epsilon : f \mapsto f(0,0,\ldots)$.

\begin{Ejercicio}
Demuestre que 
\begin{align*}
&&\Delta(h_n) &= \sum_{k+l=n} h_k\otimes h_l,\\
&& \Delta(p_n) &= p_n\otimes 1 + 1 \otimes p_n \,\,\,\,(n \ge 1)\\
&&\Delta(s_{\lambda}) &= \sum_{\mu,\nu} c_{\mu,\nu}^{\lambda} s_{\mu}\otimes s_{\nu}
\end{align*}
donde los $c_{\mu,\nu}^{\lambda}$ son los coeficientes de Littlewood-Richardson.
\end{Ejercicio}

%Se tiene que $\Sym$, junto con el producto que se hereda del producto de series formales y el coproducto obtenido a partir de la operación de suma de alfabetos, es una biálgebra.

\begin{Ejercicio}
Sobre $\Sym \otimes \Sym$ definimos un producto escalar 
\[
\langle f_1 \otimes f_2 , g_1 \otimes g_2 \rangle = \langle f_1, g_1 \rangle \langle f_2, g_2 \rangle
\]
para toda $f_1,f_2,g_1,g_2 \in \Sym$, Demuestre que
\[
\langle \Delta f, g \otimes h \rangle = \langle f, gh \rangle
\]

%Demuestre que la definición de producto escalar que acabamos de dar es equivalente  a decir que las funciones $s_{\mu}\otimes s_{\nu}$ forman una base ortonormal de $\Sym \otimes \Sym$. 
\end{Ejercicio}

%Para concluir, sabemos que $\Sym$ junto con la operación de multiplicación y su unidad tiene la estructura de álgebra graduada. Esta estructura es compatible con la estructura de coálgebra dado por $\Delta, \epsilon$. Tenemos así que $\Sym$ tiene la estructura de biálgebra. Más aún, tiene una estructura de álgebra de Hopf, donde el ant'ipoda se expresa en función de la involución  $\overline{\omega}$.
El ant'ipoda en esta bi'algebra graduada viene dada por una peque\~na variante de la involuci'on $\omega$. En efecto,
La involución 
\begin{align*}
\bar \omega : \Sym \to \Sym\\ 
\bar \omega (h_i) = (-1)^i e_i
\end{align*}
para cada $i \ge 1$, 
 es una ant'ipoda para  $\Sym$. Concluimos entonces que

\begin{Teorema}  La familia  $(\Sym, \mu, 1, \Delta, \epsilon, \bar \omega)$ es un álgebra de Hopf graduada, donde $\mu$ es la multiplicación heredada del anillo de las series formales, $1$ es su identidad, 
el coproducto viene dado por $\Delta f= f[X+Y]$,  y su counidad por $ \epsilon(f)=f(0,0,\ldots)$.

Más aún, el producto escalar de $\Sym$ es compatible con la estructura de biálgebra en el sentido de que
\begin{align*}
\langle \Delta f, g \otimes h \rangle &= \langle f, gh \rangle\\
\langle \bar \omega f, \bar \omega g \rangle &= \langle f, g \rangle\\
\langle f, 1 \rangle&= \epsilon(f)
\end{align*}

\end{Teorema}

\subsection{\textsf{La estructura de bi'algebra de $Sym$, el producto de alfabetos y los coeficientes de Kronecker.}}

Ahora pasamos a estudiar una segunda estructura de biálgebra para $Sym$ (aunque no de 'algebra de Hopf). Para esto introducimos el {\em  producto de Kronecker} 
que, como veremos,  corresponde al producto tensorial  interno entre representaciones del grupo simétrico, a trav'es de la operaci'on del producto de dos 
alfabetos.

El producto de alfabetos nos define  una segunda comultiplicación sobre $\Sym$:
\begin{align*}
\Delta^{\star} : \Sym &\to \Sym \otimes \Sym\\
\Delta^{\star} f &= f[XY].
\end{align*}
cuya counidad viene dada por
\[
\epsilon^{\star} f = f(1,0,0,\ldots)
\]
para toda $f \in \Sym.$

\begin{Ejercicio}
Demuestre que
\begin{align*}
&&\Delta^{\star} h_n = \sum_{|\lambda|=n} s_{\lambda} \otimes s_{\lambda},
&&\Delta^{\star} e_n = \sum_{|\lambda|=n} s_{\lambda} \otimes s_{\lambda'},
&&\Delta^{\star} p_n = p_n \otimes p_n.
\end{align*}
Más aún, $\epsilon^{\star} h_n = 1$, $\epsilon^{\star} e_n = \delta_{1,n}+\delta_{0,n}$ y $\epsilon^{\star} p_n = 1$.
\end{Ejercicio}

\begin{definicion}[Los coeficientes de Kronecker y el producto de alfabetos]
El desarrollo $\Delta^{\star} s_{\lambda}$ en la base de Schur:
\[
s_{\lambda}[XY]=\sum_{\mu,\nu}\kr s_{\mu}[X] s_{\nu}[Y]
\]
Nos define una nueva familia de constantes de estructuras para $Sym$,  esta vez con respecto a la comultiplicaci'on $\Delta^{\star}$. 
Los coeficientes $\kr$ los llamamos {\em coeficientes de Kronecker}. Pronto veremos que los coeficientes de Kronecker nos 
describen la multiplicidad con la que aparece la representaci'on
 irreducibe $S^{\lambda}$, en el producto tensorial $S^{\mu} \otimes S^{\nu}$. Por lo tanto, son todos enteros no--negativos.
 
 \end{definicion}
\begin{Ejercicio}
Demuestre que los coeficientes de Kronecker $\kr$ son diferentes de cero sólo si $|\lambda|=|\mu|=|\nu| $.
\end{Ejercicio}

Los coeficientes de Kronecker nos permiten definir el  {\em producto de Kronecker} entre $s_{\mu}$ y $s_{\nu}$  utilizando
para esto a la base de Schur :
\begin{equation}\label{defkron}
s_{\mu} \star s_{\nu} = \sum_{\lambda } \kr s_{\lambda}
\end{equation}
Extendiendo a $\Sym$ por linealidad, obtenemos una segunda operaci'on de multiplicaci'on para $Sym$. 
Llamamos a esta operación {\em producto de Kronecker (o producto interno)}.  

\begin{Ejercicio} Utilizando   que los elementos $s_{\mu}[X]s_{\nu}[Y]$ forman una base ortonormal de $\Sym \otimes \Sym$; demostrar que para cada $f,g,h \in Sym$
\[
\langle \Delta^{\star} f, g \otimes h \rangle = \langle f , g \star h \rangle
\]
En otras palabras, que $\Delta^{\star}$ es el adjunto del producto de Kronecker.
\end{Ejercicio}

%	\begin{Ejercicio}Sea $\lambda$ una partición de $ k$ y sea $f \in \Sym^{(k)}$.
%	Demuestre que 
%	\[
%	e_k \star f = f \star e_k = \omega f
%	\]
%	\end{Ejercicio}

\begin{Ejercicio}
Demostrar que $p_n[XY]=p_n[X]p_n[Y]$ y  $p_{\lambda}*p_{\mu}=\delta_{\lambda,\mu}z_{\lambda}p_{\lambda}$. 
Concluya que 
\[
p_{\lambda}*p_{\mu}= \langle p_{\lambda}, p_{\mu} \rangle \, p_{\lambda}.
\]
Demuestre que el producto de Kronecker es conmutativo.
\end{Ejercicio}

El siguiente ejercicio muestra que la unidad con respecto al producto de Kronecker est'a dada por una suma infinita de 
funciones homog'eneas completas, y   por
lo tanto no est'a en $Sym$ sino en su completaci'on.

\begin{Ejercicio} Demuestrar que $
h_k \star f = f \star h_k= f,$ para toda $f \in Sym^{(n)}$.
Deducir entonces que,  en la completaci'on $\widehat\Sym$ de $Sym$ la identidad con respecto al producto de Kronecker  
viene dada por $\sum_{k \ge 0} h_k$.
\end{Ejercicio}

\section{\textsf{El álgebra  de Grothendieck del grupo simétrico y del grupo lineal general.}}

En esta secci'on revisitaremos la teor'ia de representaciones del grupo sim'etrico. En la secci'on 3 ya construimos una familia completa de
representaciones irreducibles de  $\S_n$, para cada  $n$. Ahora queremos estudiar simult'aneamente estas familias de representaciones
irreducibles desde un punto de vista algebraico.

Sea $\grn$ el grupo de Grothendieck del grupo simétrico $\S_n$. Esto es, el grupo abeliano generado por las clases de equivalencia de las representaciones irreducibles del grupo simétrico con la operación de  suma directa.
Definimos 
\[
\gr = \{ \grn : n \ge 1 \}
\]
junto con la inclusión canónica: 	
\[
\rho_{n,m} : \grn \times R(\S_m) \hookrightarrow R(\S_{n+m}).
\]
que considera un par de  permutaciones $(\pi,\sigma)$ en $\S_n\times\S_m$ como una permutaci'on de $\S_{m+m}$ aplicando 
 $\pi$ a los n'umeros $1, \cdots, n$ y $\sigma$ a 
$n+1, \cdots, m$ de la manera can'onica : $\sigma(n+i):=\sigma(i)$.

Sobre $\gr$ se puede definir una estructura de álgebra de Hopf y una segunda estructura de biálgebra que, como veremos, corresponden a  las estructuras que acabamos de estudiar sobre $\Sym$.

\begin{definicion}[La estructura de álgebra de Hopf sobre $\gr$]

Sean $X$ e $Y$ dos  $\S_n$--módulos. Definimos {\em el producto entre $X$ e $Y$}, que denotamos por $X \circ Y$, a través del procedimiento de inducción de representaciones:
\begin{align}\label{producto}
X \circ Y := \rho_{n,m} (X \otimes Y) \big\uparrow_{\S_n \times \S_m}^{\S_{n+m}}
\end{align}

%\todo{Los coeficientes de LR}

{\em La unidad} de esta álgebra de Hopf viene dada por la unidad del álgebra $R(\S)$; es decir, la representación trivial de $\S_0$.

{\em El coproducto} viene dado por la suma de las restricciones de  $R(\S_n) $ a $R(\S_i) \otimes R(\S_j)$, para $i+j=n$  enteros no negativos.
\begin{align*}
\Delta : R(\S_n) &\to \sum_{i+j=n}  R(\S_i) \otimes R(\S_j)\\
\Delta(X) &=  \sum_{i + j = n}   X  \big\downarrow^{\S_{i+j} }_{\S_i \times \S_j}
\end{align*}

{\em La counidad} está definida  por la  proyección sobre $R(\S_0)$ y nos da la multiplicidad de la representación trivial de $\S_0$.

Nótese que no es trivial que la comultiplicación sea un morfismo de álgebras. Esto es una consecuencia del teorema de Mackey. Para una exposición completa de estos resultados ver el libro de Zelevinski \cite{Zelevinski}.

El ant'ipoda corresponde (salvo por un signo) a tomar el producto tensorial con la representación alternante.

%\todo{hablar sobre el teorema de Mackey}

\end{definicion}

 \begin{definicion}[La estructura de biálgebra  sobre $\gr$ con respecto al producto de Kronecker]
 Sean $X$ e $Y$ dos $\S_n$-módulos. El producto de Kronecker (también llamado producto tensorial interno) se define  como el $\S_n$--módulo $X \otimes Y$ junto con la acción diagonal del grupo simétrico. Esto es, para cada $\pi \in \S_n$
 \[
 \pi ( X \otimes Y) = \pi X \otimes \pi Y 
 \] 
 Denotamos el $\S_n$--módulo resultante por $X \star Y$. Nótese que si $\chi$ es el car\'acter de $X$ y $\phi$ el car\'acter de $Y$, entonces el car\'acter del producto de Kronecker viene dado por $(\chi\phi) (\pi) = \chi(\pi) \phi(\pi)$.

%	 \begin{Ejercicio}
%	 Demuestre que el producto de Kronecker es asociativo, conmutativo y que su unidad viene dada por 
%	 FALTA su imagen en $h_n$.
%	 \end{Ejercicio}

Los  {\em coeficientes de Kronecker} tambi'en pueden ser definidos como las multiplicidades de las representaciones irreducibles en el producto tensorial de dos representaciones irreducibles del grupo simétrico. Luego veremos que esta definici'on es equivalente a la que dimos en la secci'on anterior en t'ermino de
operaciones entre alfabetos.

Si $\mu$ y $\nu$ son particiones de $n$, se definen los coeficientes de Kronecker $\kr$ como
\[
\chi^{\mu}\chi^{\nu}=\sum_{\lambda \vdash n} \kr \chi^{\lambda}.
\]

%\todo{Hablar de los caracteres}
Equivalentemente, la siguiente ecuación define los coeficientes de Kronecker y demuestra que son simétricos en $\mu, \nu $ y $\lambda$:
\[
\kr=\langle \chi^{\lambda}, \chi^{\mu}\chi^{\nu} \rangle_{\S_n}= \frac{1}{n!}
\sum_{\sigma \in \S_n}  \chi^{\lambda}(\sigma) \chi^{\mu}(\sigma) \chi^{\nu}(\sigma),
\]
donde $\langle \,\,\, , \, \rangle_{\S_n}$ es el producto escalar en el espacio generado por los caracteres de $\S_n$, sus funciones de clase.\footnote{La barra de conjugaci'on desaparece en este producto porque todas las representaciones irreducibles de $\S_n$  pueden ser construidas sobre $\mathbb{Q}$.}

 Extendemos la definición del producto de Kronecker a $R(\S)$ diciendo que el producto de Kronecker entre una  representación $\S_n$ y otra  de $\S_m$ es cero si $n \ne m$.  
 \begin{Ejercicio}
 Demuestre que la unidad de $\Sym^{(n)}$, con respecto al producto de Kronecker, está dada por la representación trivial de $\S_n$. Describa la counidad y la comultiplicación correspondiente.
 \end{Ejercicio}

 \end{definicion}
 
Es importante darse cuenta que si $X$ es un $\S_n$--m'odulo y $Y$ es un $\S_m$--m'odulo, entonces $X \circ Y$ es un $\S_{m+n}$--m'odulo
 y $X * Y$ es un $\S_n$-m'odulo.

\section{\textsf{La aplicación de Frobenius}}

%%%%%%%%%%%%%%%%%%%%%%%%%%%%%%%%%%%%%%%%%%%%%%%%%
%%%%%%%%%%%%%%%%%%%%%%%%%%%%%%%%%%%%%%%%%%%%%%%%%
%%%%%%%%%%%%%%%%%%%%%%%%%%%%%%%%%%%%%%%%%%%%%%%%%
%%%%%%%%%%%%%%%%%%%%%%%%%%%%%%%%%%%%%%%%%%%%%%%%%
%%%%%%%%%%%%%%%%%%%%%%%%%%%%%%%%%%%%%%%%%%%%%%%%%
%%%%%%%%%%%%%%%%%%%%%%%%%%%%%%%%%%%%%%%%%%%%%%%%%
%%%%%%%%%%%%%%%%%%%%%%%%%%%%%%%%%%%%%%%%%%%%%%%%%
%%%%%%%%%%%%%%%%%%%%%%%%%%%%%%%%%%%%%%%%%%%%%%%%%

Ahora nos proponemos explicar la relación existente entre la teoría de representaciones del grupo simétrico y la funciones simétricas $\Sym$. 
Las funciones sim'etricas juegan el rol de funciones generatrices para los caracteres de las representaciones del grupo sim'etrico.  
%Para despertar la curiosidad del lector, le planteamos el siguiente ejercicio:

%\begin{Ejercicio}
%Verificar las siguientes igualdades entre funciones sim'etricas, y compararlas con  los coeficientes que aparecen en el Ejercicio \ref{caracteres}:
%\begin{align*}
%p_3           &= m_{3}\\
%p_{2, 1}    & = m_3+m_{2,1}\\
%p_{1,1,1}  &= m_3 + 3 m_{2,1} + 6 m_{1^3}.
%\end{align*}
%\end{Ejercicio}

%Parece ser que el desarrollo de las funciones sim'etricas series de potencias en la base monomial nos proporcionan los  caracteres de las representaciones irreducibles del grupo simétrico. El siguiente teorema nos asegura que esto es cierto:

%\begin{Teorema}\label{p2m}
% Sean $\lambda$ y $\mu$ particiones de $n$.
%Sea $\chi^{\mu}({\lambda})$ el car\'acter de la representaci'on irreducible $S^{\mu}$ evaluado en la clase que corresponde a $\lambda$. Entonces
%\[
%p_{\lambda} = \sum_{\mu \unrhd \lambda} \chi^{\mu}({\lambda}) m_{\mu}
%\]
%\end{Teorema}
%
%Por otra parte, el desarrollo de la base de Schur en la base de las sumas de potencias es particularmente relevante en nuestro estudio.

\begin{Teorema}[Frobenius] \label{s2p}
Si $\lambda \vdash n$, se tiene entonces que
\[
s_{\lambda} = \frac{1}{n!} \sum_{\pi \in \S_n} \chi^{\lambda}(\pi) p_{\pi}=\sum_{\mu \vdash n} \chi^{\lambda}(\mu) \frac{p_{\mu}}{z_{\mu}}.
\]
\end{Teorema}
Este resultado central nos describe elegantemente la relaci'on que existe entre la teor'ia de las funciones sim'etricas, y la teor'ia de 
representaciones del grupo sim'etrico,   \cite{F, Sagan, EC1}. Es este resultado el que nos permite definir  la aplicaci'on caracter'istica de Frobenius.

Sea $CF^k$ el conjunto de todas las funciones de clase $f : \S_k \to \mathbb{Q}$ (funciones constantes en las clases de conjugaci'on de 
$\S_k$.) En $CF^k$ existe un producto escalar natural definido por
\[
\langle f, g\rangle = \frac{1}{n!} \sum_{\pi \in \S_k} f(\pi)g(\pi)
\]
Escribimos $\langle X, Y \rangle$ cuando $X$ y $Y$ son representaciones de $\S_k$ con car'acteres $\phi$ y $\chi$, respectivamente. 

Estamos interesados en estudiar la siguiente transformaci'on lineal :
\begin{definicion}[La caracter'istica de Frobenius]  
Sea $f$ es una funci'on de clase de grado $k$. Definimos la aplicaci'on
\begin{align*}
&\ch^{(k)} : CF^k  \to Sym^{(k)}\\
&\ch^{(k)}(f) = \sum_{\mu \vdash k}f({\mu}) \frac {p_{\mu}}{ z_{\mu}}
\end{align*}
donde $f(\mu)$ denota $f(\pi)$ para cualquier permutaci'on $\pi$ de tipo $\mu$.
%donde  $\chi(\mu)$ denota el valor del car\'acter $\chi$ del $\S_k$-módulo $X$ en la clase de conjugación $\mu$. La extensión natural de $\ch^{(k)}$ sobre $R(\S)$ la denotamos por $\ch.$

\end{definicion}

El teorema de Frobenius nos dice entonces que  
\[
\ch(\chi^{\lambda})=s_{\lambda}.
\]

Este importante resultado nos permite  definir las funciones de Schur a través del importante rol que juegan en la teoría de representaciones: {\em Las funciones de Schur son las imágenes, bajo la aplicación de Frobenius, de las representaciones  irreducibles.}

El resultado anterior implica que $\ch$ es un isomorfismo de espacios vectoriales. La aplicaci'on de Frobenius env'ia una base ortonormal de $CF^k$, los 
caracteres de las representaciones irreducibles de $\S_n$, en la base ortonormal de $\Sym^{(k)}$ : las funciones de Schur. Concluimos entonces que
la aplicaci'on de Frobenius es una isometr'ia.

Podemos utilizar la aplicaci'on de Frobenius para obtener informaci'on sobre las funciones sim'etricas a partir de nuestros conocimientos de la
teor'ia de representaciones del grupo sim'etrico. Por ejemplo, como es inmediato calcular los caracteres de las representaciones triviales y alternadas
podemos concluir que

\begin{Ejercicio}[Frobenius, y las representaciones trivial y alternante.] Utilizar que conocemos los caracteres de las representaciones trivial 
y alternante para demostrar que
\begin{align*}
h_{n} &= \sum_{\mu\vdash n} s_{\mu},\\
e_{n} &= \sum_{\mu\vdash n} (-1)^{signo(\mu)}s_{\mu}.
\end{align*}
donde el signo de una partici'on $\mu$ se define como el signo de cualquier permutaci'on de tipo $\mu$.
\end{Ejercicio}

%Veamos ahora algunas propiedades de la aplicación de Frobenius.  Recordemos que hemos definido dos estructuras sobre $\Sym$, una de álgebra de Hopf con respecto al producto de series, $(\Sym, \cdot)$, donde $\cdot$ denota al producto de series,  y una segunda estructura de biálgebra con respecto al producto de Kronecker, $(\Sym, \star)$.

% Similarmente, hemos definidos dos estructuras sobre $R(\S)$. La primera, una estructura de álgebra de Hopf con respecto al producto $\circ$ definido en función de la operación de inducción de representaciones,  $(R(\S), \circ)$, y una segunda estructura de biágebra definida en función del producto de Kronecker $(R(\S), \star).$

\begin{Teorema}[Frobenius]
La aplicación de Frobenius  es multiplicativa. Más precisamente:
\begin{align*}
\ch (X \circ Y) &= \ch(X) \phantom{\star} \ch(Y),\\
\ch (X \star Y) &= \ch(X) \star \ch(Y).
\end{align*}
\end{Teorema}
La primera afirmación es una consecuencia de la fórmula de inducción de Frobenius. La segunda se sigue inmediatamente de la definición del producto de Kronecker. 

\begin{Ejercicio}
Describa el efecto de la aplicación de Frobenius sobre la unidad, la comultiplicación , la counidad y el ant'ipoda del álgebra de Hopf de $R(\S_k)$ .
\end{Ejercicio}

\begin{Teorema} La aplicación de Frobenius 
\[
\ch : (R(\S), \circ)  \to (\Sym, \cdot) 
\]
es un isomorfismo de álgebras de Hopf. Y
\[
\ch : (R(\S), \star)  \to (\Sym, \star) 
\]
es  un isomorfismo de biálgebras. 
\end{Teorema}

%\begin{Ejercicio}[La aplicación de Frobenius es una isometría]
%Demuestre  que $\ch^{(k)}$ preserva el producto escalar.  \end{Ejercicio}
\subsection{\textsf{Los caracteres de las representaciones irreducibles de grupo lineal general}}

Sea $V$ un espacio vectorial de dimensi'on $m$. Escogiendo una base para $V$ identificamos $GL(V)$ con $GL_m(\mathbb{C})$. Sea $H$ el subgrupo de las
matrices diagonales de $GL_m(\mathbb{C})$, y sea $diag(X)$ la matriz diagonal cuyas entradas en la
diagonal principal son $(x_1, \cdots, x_n)$.

Un vector $v$ en una representacion $V$ se llama vector de peso, con peso $\alpha=(\alpha_1, \alpha_2, \cdots, \alpha_m)$, con coordenadas enteras, si
\[
x \cdot v = x_1^{\alpha_1} x_2^{\alpha_2}\cdots x_m^{\alpha_m} v,  \text{ para todo $x$ in $H$}
\]
Se tiene entonces que cualquier representaci'on $V$ es una suma directa de espacios de peso $V= \oplus V_{\alpha}$, donde $V_{\alpha} = \{v \in V : x \cdot v = x_1^{\alpha_1} x_2^{\alpha_2}\cdots x_m^{\alpha_m} \, v, \text{para todo $x \in H$} \}$.

Una representaci'on (finito-dimensional y holomorfa) $V$ de $GL_m(\mathbb{C})$ es irreducible si y s'olo si tiene un \'unico vector de peso $\alpha$ en la descomposici'on
anterior. M'as a'un, dos representaciones son isomorfas si  y s'olo si tienen el mismo vector de peso. Hemos visto que las representaciones polinomiales de $GL(V)$ vienen
indexadas por particiones, y se pueden construir a partir de las representaciones irreducibles del grupo sim'etrico.

El car'acter de una representaci'on (finito dimensional y holom'orfica) 
$W$ de $GL(V)$,  $Char_W=Char_W(x_1,\cdots,x_m)=\chi_W$ se define como la  traza de $diag(x)$ sobre $W$.

Si descomponemos a $V= \oplus V_{\alpha}$ tenemos entonces que 
\[
\chi_V(x)=\sum_{\alpha}\dim(V_\alpha)x^\alpha=\sum_\alpha \dim(V_\alpha) x_1^{\alpha_1} x_2^{\alpha_2}\cdots x_m^{\alpha_m}
\]
Y en particular, para $V^{\lambda}$ tenemos un vector de peso para cada tableaux $T$ con entradas en $[m]$, de manera que
\[
Char(V^{\lambda})=\sum X^T = s_{\lambda}(x_1, x_2, \cdots, x_m)
\]
es el polinomio de Schur indexado por $\lambda$.

En general, tenemos que
\begin{align*}
%Char(W^\lambda)&= \sum_T x^T = s_{\lambda}(x_1, x_2, \cdots, x_n)\\
Char(W\oplus W') &=Char(W) \oplus Char(W')\\
Char(W\otimes W') &= Char(W) \, Char(W')
\end{align*}

\subsection{\textsf{El algebra de Grothedieck del grupo lineal general y el fen'omeno de dualidad de Schur-Weyl.}} 
Definimos  el anillo de representaciones de $
R(GL_m)$ como el grupo abeliano generado por las clases de isomorf'ia de las representaciones polinomiales de $GL_m$, junto con la operaci'on de la suma
directa. 

Veamos que $R(GL_m)$ hereda  de la operaci'on de inducci'on entre representaciones
del grupo sim'etrico una estructura de 'algebra de Hopf con respecto al producto tensorial de representaciones.
 
 Sea $X$ es una representaci'on de $\S_n$ y $Y$ una representaci'on de $\S_m$, hemos visto que 
 $X \circ Y $ es una representaci'on de $\S_{n+m}$ . Es sencillo demostrar que 
 \[
 \mathbb{V}(X\circ Y) \cong \mathbb{V}(X) \otimes \mathbb{V}(Y)
 \]
 Por lo que la operaci'on de inducci'on de representaciones del grupo sim'etrico corresponde a tomar el producto tensorial de
 representaciones de $GL(E).$
M'as a'un, el problema de descomponer $S^{\lambda}\circ S^{\mu}$ como suma de representaciones irreducibles de $\S_{n+m}$, es equivalente
al problema de descomponer $V^{\lambda} \otimes V^{\mu}$. En ambos casos los coeficientes de estructura que aparecen son los
coeficientes de Littlewood-Richardson.

 Esto evidencia una profunda relaci'on entre las teor'ia de representaciones del grupo sim'etrico y del grupo lineal general, la llamada 
 {\em dualidad de Schur-Weyl.}  Sea $V$ un espacio vectorial de dimensi'on $m$. Recordemos que hemos definido sobre $V^{\otimes m}$
 una estructura de $(GL(V), \S_n)$--m'odulo, donde el grupo $GL(V)$ act'ua diagonalmente, y $\S_n$ act'ua por la derecha.
 Estas dos acciones conmutan entre s'i, y nos proporcionan la descomposici'on.
 \[
 V^{\otimes m} = \sum_{\substack{\lambda\vdash n\\\ell(\lambda)\le m}} S^{\lambda} \otimes V^{\lambda}
 \]
 Esta descomposici'on nos proporciona una manera de estudiar las representaciones del grupo sim'etrico a trav'es de aquellas del 
 grupo lineal general, y viceversa. Por ejemplo, supongamos que queremos interpretar el producto de Kronecker, entre representaciones del grupo
 sim'etrico, en el lenguaje de la teor'ia de representaciones del grupo lineal general.  Sea $V$ un espacio vectorial de dimensi'on $n$ y $W$ un
 espacio vectorial de dimensi'on $m$. Para hacer nuestra escritura más transparente, denotamos por $V^{\lambda}_m$ a la representaci'on
 irreducible de $GL_m$ indexada por $\lambda$.
 
 Tenemos entonces que $\dim(V\otimes W)=mn$ y
 \begin{align*}
 (V\otimes W)^{\otimes k} = \sum_{\lambda\vdash k, \ell(\lambda)\le nm} S^{\lambda} \otimes V^{\lambda} 
 \end{align*}
 donde  $V_{nm}^{\lambda}$ es la representaci'on irreducible de $V\otimes W$ indexada por $\lambda$. Esta representaci'on no es
 irreducible. Supongamos que se descompone de la forma:
 \[
 V^{\lambda}_{n,m} = \sum _{\substack{\mu\vdash m, \ell(\mu)\le n\\ \nu \vdash m, \ell(\nu)\le m }}  (V^{\mu}_n\otimes V_m^{\nu})^{g_{\mu,\nu}^{\lambda}}.
 \]
 Aplicando otra vez la dualidad de Schur-Weyl, a $ (V\otimes W)^{\otimes k}$, pero interpretando esta expresi'on como $V^{\otimes k} \otimes W^{\otimes k}$,
 obtenemos
  \[
 V^{\lambda}_{n,m} = \big(\sum _{\mu\vdash k,\, \ell(\mu) \,\le n}  V^{\mu}_n\otimes S^{\mu}\big)\otimes  \big(\sum _{\nu\vdash k, \ell(\nu)\le m}  V^{\nu}_m\otimes S^{\nu}\big)
 \]
 Todas estas ecuaciones son ciertas, m'odulo isomorfismo, as'i que podemos utilizar  la propiedad conmutativa, junto con la descomposici'on del
 producto tensorial de representaciones irreducibles del grupo sim'etrico para concluir que los coeficientes $g_{\mu,\nu}^{\lambda}$ son precisamente
 los coeficientes de Kronecker.
 \begin{definicion}[Los coeficientes de Kronecker y el grupo sim'etrico] Los coeficientes de Kronecker son las constantes de estructura que
 aparecen al descomponer el producto tensorial de representaciones irreducibles del grupo sim'etrico.
 \[
S^{\mu} \otimes S^{\nu} = \sum_{\lambda\vdash k} (S^{\lambda})^{\otimes \gamma_{\mu,\nu}^{\lambda}}
 \]
 N'otese que esta ecuaci'on es cierta para cada valor de $nm$. Por lo tanto podemos omitir la condici'on  $\ell(\lambda)\le nm$.
 \end{definicion}
 
% En el lenguage de la teor'ia de invariantes
% \[
%\gamma_{\mu,\nu}^{\lambda} = \dim((V^*)^{\lambda}\otimes V^{\mu}\otimes V^{\nu})^{GL(V)\times GL(W)}=\dim(S^{\lambda}\otimes S^{\mu}\otimes S^{\nu})^{\S_k}
%\]

 Finalmente, podemos utilizar la dualidad de Schur-Weyl para interpretar los coeficientes de Kronecker en el lenguage de la
 teor'ia de representaciones del grupo sim'etrico. Los coeficientes $\gamma_{\mu,\nu}^{\lambda}$ describen
  la descomposici'on en irreducibles de la representacion $V^{\lambda}$
 de $GL(V \otimes W)$ como suma de representaciones irreducibles de $GL(V)\times GL(W)$. Finalmente, utilizando la f'ormula de reciprocidad
 de Frobenius podemos describir estos coeficientes en funci'on de la operaci'on de inducci'on.

  \section*{\textsf{Agradecimientos}}

Agradecemos a Emmanuel Briand por su apoyo y por varias discusiones matem\'aticas de gran utilidad.

  \end{document}